\documentclass[12pt]{amsart}
\usepackage{amssymb, latexsym, graphicx}
\usepackage{graphicx}
\usepackage{amsmath}
\usepackage{amsthm}
\usepackage{lscape}
\usepackage{booktabs}
\usepackage[color]{}
\usepackage{fancyhdr}

\def\ep{\varepsilon}
\def\f{\frac}

\def\p{\partial}

\numberwithin{equation}{section}

\theoremstyle{remark}
\newtheorem{remark}{Remark}
\numberwithin{remark}{section}

\begin{document}
\title[Treatment of Incompatibilities]{Treatment of incompatible initial and boundary data for parabolic equations in higher dimension}

\author{Qingshan Chen}
\address[QC]{Department of Scientific Computing\\
         The Florida State University\\
	 Tallahassee, FL 32306}
\email{qchen3@fsu.edu}
\author{Zhen Qin}
\address[ZQ]{The Institute for Scientific Computing and Applied Mathematics\\
         Indiana University\\
	 Bloomington, IN 47405}
\email{qinz@indiana.edu}
\author{Roger Temam}
\address[RT]{The Institute for Scientific Computing and Applied Mathematics\\
         Indiana University\\
	 Bloomington, IN 47405}
\email{temam@indiana.edu}

\date{\today}
\maketitle
\thispagestyle{fancy}
\renewcommand{\headrulewidth}{0pt}
\fancyhead{ }
\fancyhead[LO]{\textsf{Article to appear in {\em Math.~Comp.}}} 





\begin{abstract}
A new method is proposed to improve the numerical simulation of time dependent problems
when the initial and boundary data are not compatible.  Unlike earlier methods limited
to space dimension one, this method can be used for any space dimension.
When both methods are applicable (in space dimension one), the improvements in precision
are comparable, but the method proposed here is not restricted by dimension.
\end{abstract}
\section{Introduction}\label{s1}
When performing large scale numerical simulations for evolutionary problems, we use most often initial and boundary conditions provided by approximations, by other simulations, or by experimental measurements. These data may not satisfy certain compatibility conditions verified by the solutions; thus various modifications deemed non essential are made on the data to overcome these difficulties. Such issues are extensively addressed in the literature; see for instance in geophysical fluid mechanics \cite{BF99} or \cite{Tren} which contains many allusions to this difficulty; in classical fluid mechanics, see e.g. \cite{ESG82,Gre91,GS87,Hey80,HR83}; see also \cite{Bie05} in chemistry and \cite{TW74} in a general mathematical context.\\

We want to address here a less known difficulty of "mathematical" nature which, the specialists believe, will become very important as we move to high resolution methods thanks to the increase of computing power and memory capacity of the computers. A very simple example of such a difficulty appears when solving in space dimension one on $(0,1)$, the heat equation $u_{t}-u_{xx}=0$ with boundary conditions $u(0,t)=u(1,t)=0$ and initial condition $u(x,0)=1$. The solution exists and is unique (for $t>0$) and the analytic expressions of $u$ are provided in the literature (see e.g. \cite{Ca84}). This problem is simple enough that it can be solved satisfactorily by numerical methods, but the solution does display singularities in the corner $x=0,t=0$ and $x=1,t=0$. For this problem and for general parabolic equations, it is known from semi-group theory \cite{HE81,Pa83} or by using the analyticity in time of the solutions (see \cite{FoTe79}) that certain norms of $\partial u /\partial t$ grow as a power of $1/t$ when $t\rightarrow 0$. It is believed that such singularities will affect large scale computations as our demand for better results increases. In fact it has been observed by some authors that, when using spectral methods for the space discretization, the spectral accuracy is lost if nothing is made to address this singularity and a series of works resulted from this observation; see e.g. \cite{BoFl99,BF99,FlFo03,FlFo04,FlSw02}, see also \cite{CQT11} for a nonlinear equation.\\

The mathematical difficulty studied in detail in e.g. \cite{LaSoUr68,La52,La54,RaMa74,Sm80,Te82} is the following; even if the initial and boundary data of an evolution problem are given $\mathcal{C}^{\infty}$, the solution may not be $\mathcal{C}^{\infty}$ near $t=0$. In fact, $k$ compatibility conditions between the data are needed for the solution to be $\mathcal{C}^{k}$ near $t=0$ and hence an infinite number of compatibility conditions are needed for the solution to be $\mathcal{C}^{\infty}$. Furthermore the initial and boundary conditions that are compatible form a relatively small set (in an informal sense), so that most numerical simulations are done with data which are not compatible, generating a loss of accuracy near $t=0$ if nothing is done. In the works mentioned above, methods have been proposed to address the first or the first two incompatibilities. It is believed that dealing with one or two incompatibilities substantially improve the quality of the simulation and, in any case, dealing with more incompatibility conditions may become impractical. However a strict limitation of these works is that the proposed methods only apply to space dimension one and, to the best of our knowledge, there is (there was) no method available in dimension two or larger to address this difficulty.\\

As we said, past works, on the computational side, have been devoted to space dimension one. This problem has been addressed in a series of articles by Flyer, Boyd, Fornberg and Swarztrauber \cite{BoFl99, FlSw02, FlFo03, FlFo04} who proposed a number of remedies in space dimension one for linear equations.  Nonlinear equations in space dimension one were considered in \cite{CQT11}.\\

In these articles, the authors introduce a correction term in the linear and nonlinear cases
by setting
\begin{equation}\label{e1.1}
u=v+S,
\end{equation}
where $S$ absorbs the incompatibilities between the initial and boundary data
up to a certain order.
Now, free of  incompatibilities of lower orders (the most severe ones), $v$
is computed by an appropriate numerical procedure, such as finite differences, the Galerkin finite
element method, spectral or pseudo--spectral methods. As a final step,
the original solution $u$ is recovered through \eqref{e1.1}. This remedy
procedure effectively reduces the errors at the spatio--temporal corners
during the short initial transient period.\\

For dimensions higher than one, the construction of $S$, to correct for singularities generated at $t=0$ by incompatible data, remains an open problem. A method to overcome this difficulty is proposed, analyzed and tested in this article.\\

In this article, we intend to study, from both a theoretical and
numerical point of view, the incompatibility issue for the
multi-dimensional time-dependent linear parabolic equation:
\begin{equation}\left\{\label{e1.2}
\begin{aligned}
&u_{t}-\nu \triangle u= f,\qquad x\in\Omega\subset R^{d}, \quad t\in R^{+},\\
&u|_{t=0}=u_{0},\\
&u|_{\partial\Omega}=g.
\end{aligned}\right.
\end{equation}

We believe that our method applies to more general parabolic equation, but we restrict ourselves to equation (\ref{e1.2}) in this article devoted to feasibility.\\

The method that we propose is based on the concept of penalty. We replace in (\ref{e1.2}) the boundary value $u|_{\partial\Omega}=g$ by $u|_{\partial\Omega}=k^{\varepsilon}$. This boundary value $k^{\varepsilon}$ which depends on a parameter $\varepsilon>0$ is such that $k^{\varepsilon}|_{t=0}=u_{0}|_{\partial\Omega}$ (see equations (\ref{e2.1}), (\ref{e2.2}) below), so that the first incompatibility has disappeared. Now $k^{\varepsilon}$, through a penalty procedure with parameter $\varepsilon$, is forced to rapidly vary from $u_{0}|_{\partial\Omega}$ at $t=0$ to the desired value, namely $g$. This is achieved through equation (\ref{e2.2}); initially $k^{\varepsilon}_{t}$ is large, but it becomes rapidly of order 1, and then, by the first equation (\ref{e2.2}), $k^{\varepsilon}-g$ is of order $\varepsilon$. It is easy of course to integrate equation (\ref{e2.2}) although an explicit solution is not available in the general case where $g$ depends on time. The concept of penalty has been introduced in the mathematical literature by R. Courant \cite{Co43}; it has been adapted to evolution problems by J. L. Lions in \cite{Li69}, a reference which contains many evolution equations similar to (\ref{e2.2}) (Chapter 3, Sections 5 to 8); it is widely used in optimization\
\footnote{A search on Google with the words "optimization, penalty" produced 3,350,000 entries.}; see also \cite{Te01} (Chapter 1, Section 6). In this work, we firstly present our approach in details and study it theoretically to prove the
strong convergence of the method. Then we implement it numerically on a number of examples. Because the penalty method does
not depend on the properties of $\Omega$, we believe
that this method can be applied to many systems with many different domains $\Omega$. The question that remains is the choice of $\varepsilon>0$ small. In optimization theory, the choice of $\varepsilon$ is usually made by trial and error and is not a major issue. It does not follow the "intuitive" idea that the error becomes smaller as $\varepsilon$ becomes smaller because of many other contingent errors such as round-off and descretization errors. In general the error becomes "optimal" for some value of $\varepsilon$ and the method gives less good results for smaller or larger values of $\varepsilon$. In our case (see Fig. \ref{f5} and \ref{f4.1}),  at the initial steps, the error decreases sharply as $\varepsilon$ increases and remains close to $0$, then it becomes stable flat. At the final steps, the error increases almost linearly as $\varepsilon$ increases. With $\varepsilon$ at about $0.1$, the initial error is minimized while the error at the final step is well controlled.
In a short time period $\varepsilon=0.5$ gives us smaller errors and again after a short time period $\varepsilon=0.1$ gives us a smaller errors. In general the choice of $\varepsilon$ really depends on our goals of the computation.\\

This article is organized as follows. In Section 2 we present the method and establish various approximation results.
Then, in Section 3 we present numerical results showing the efficiency of
the method and comparing it to earlier methods. In Section 4 we present
some conclusions and perspective of future developments.

\section{penalty method}\label{s2}
\subsection{Perturbed problem (and the statement of the main result)}
We consider the system (\ref{e1.2}), where $\nu>0$. If
$u_{0}|_{\partial \Omega}\neq g(0)$, then we face an incompatibility problem, in which
case we consider a new system instead, namely, for $\varepsilon>0$ fixed,
\begin{equation}\left\{\label{e2.1}
\begin{aligned}
&u^{\varepsilon}_{t}-\nu\triangle u^{\varepsilon}=f, \qquad x\in\Omega\subset R^{d}, \quad t\in R^{+},\\
&u^{\varepsilon}|_{t=0}=u_{0},\\
&u^{\varepsilon}|_{\partial\Omega}=k^{\varepsilon}.
\end{aligned}\right.
\end{equation}
\begin{equation}\left\{\label{e2.2}
\begin{aligned}
&k^{\varepsilon}_{t}+\frac{1}{\varepsilon}(k^{\varepsilon}-g)=0,\qquad t\in R^{+},\quad\\
&k^{\varepsilon}(0)=u_{0}|_{\partial\Omega}.
\end{aligned}\right.
\end{equation}

In this article, $|\cdot|$ is the $L^{2}(\Omega)$ norm, and $\parallel \cdot\parallel = |\nabla\cdot|$ is the $H^{1}_{0}(\Omega)$ norm; for other norms, we will use the subscript notation.\\

The system (\ref{e2.1})-(\ref{e2.2}) is actually decoupled and (\ref{e2.2}) is
just an Ordinary Differential Equation with $x\in\partial\Omega$ as a parameter. As we see below, if we are given $g$, $g'=\dfrac{\partial g}{\partial t} \in
L^{2}(0,T;H^{\frac{1}{2}}(\Gamma))$, then we have the existence and
uniqueness of $k^{\varepsilon}$ in $L^{2}(0,T;H^{\frac{1}{2}}(\Gamma))$ and furthermore, by the effect of the penalty term, $(k^\varepsilon - g)/\varepsilon,k^{\varepsilon}$ converges to $g$ in
suitable spaces as $\varepsilon \rightarrow 0.$  Equations
(\ref{e2.1}) is a heat equation with non--homogeneous boundary conditions, and we have the
existence and uniqueness of a solution if the data are sufficiently regular. Then we have the following theorem.\\

\vspace{0.1cm} \textbf{Theorem 2.1}.\label{t1}\textit{\quad Assume that we are given $g\in
L^{\infty}(0,T;H^{\frac{1}{2}}(\Gamma))\enspace(\Gamma=\partial\Omega),$ with
$g_{t}\in
L^{2}(0,T;H^{\frac{1}{2}}(\Gamma))$, and $u_{0}\in H^{1}(\Omega).   $   Then (\ref{e1.2}) has a unique solution $u\in L^{2}(0,T;H^{1}_{0}(\Omega))$ $ \cap$ $  \mathcal{C}([0,T];L^{2}(\Omega))$, and for each $\varepsilon>0$, (\ref{e2.1})-(\ref{e2.2}) has a unique solution $u^{\varepsilon}\in L^{2}(0,T;H^{1}(\Omega)) \cap$ $\mathcal{C}([0,T];$ $L^{2}(\Omega))$, $k^{\varepsilon}\in L^{2}(0,T;H^{\frac{1}{2}}(\Gamma))$. Furthermore, as $\varepsilon\rightarrow 0$,}

\begin{equation}\label{e2.3}
\begin{aligned}
u^{\varepsilon}\rightarrow u\enspace in\enspace &L^{2}(0,T;H^{-1}(\Omega))\enspace strongly, and\enspace in \\
&\mathcal{C}([t_{0},T];H^{-2}(\Omega))\enspace strongly,\enspace \forall t_{0}>0.
\end{aligned}
\end{equation}
\begin{remark}\label{r2.1}
We do not prove a strong convergence of $u^{\varepsilon}$ to $u$ on all of $[0,T]$ in
the $L^{\infty}$ sense, and we do not expect such a convergence to occur since $u$
has a singularity at $t=0$. Alternatively one could capture the singularity of
$u$ near $t=0$ by using the methods of singular perturbation theory as in,
e.g.~Jung-Temam \cite{JuTe05}, which we do briefly in
Section \ref{s2.2}, and will also be studied elsewhere.
\end{remark}

Before we prove Theorem 2.1, we will first prove the following lemma.\\

\textbf{Lemma 2.1}.\label{l2.1} \textit{If $g\in L^{\infty}(0,T;H^{\frac{1}{2}}(\Gamma))$ and  $g'\in L^{2}(0,T;H^{\frac{1}{2}}(\Gamma))$, then there exists a unique $k^{\varepsilon}$ in  $L^{2}(0,T;H^{\frac{1}{2}}(\Gamma))$ satisfying (\ref{e2.2}), and as $\varepsilon \rightarrow 0$, $k^{\varepsilon}\rightarrow g $ in $L^{2}(0,T;H^{\frac{1}{2}}(\Gamma))$ strongly. Furthermore, as $\varepsilon\rightarrow 0$, $\displaystyle \int^{t}_{0}k^{\varepsilon}(s)ds\rightarrow \int^{t}_{s}g(s)ds$  in $L^{2}(0,T;H^{\frac{1}{2}}(\Gamma))$ strongly.}\\

\begin{proof}
We explicitly solve the ODE system (\ref{e2.2}), and we obtain the solution $k^{\varepsilon}\in L^{2}(0,T;H^{\frac{1}{2}}(\Gamma))$:
\begin{equation}
k^{\varepsilon}(t)=e^{-\frac{t}{\varepsilon}}k^{\varepsilon}(0)+\int^{t}_{0}\frac{1}{\varepsilon}g(s)e^{\frac{s-t}{\varepsilon}}ds.
\end{equation}

Then we rewrite $(\ref{e2.2})_{1}$ in the form\\
\begin{equation}\label{e2.5}
(k^{\varepsilon}-g)_{t}+\frac{1}{\varepsilon}(k^{\varepsilon}-g)=-g_{t}.
\end{equation}

Taking the scalar product of (\ref{e2.5}) with $k^{\varepsilon}-g$ in $H^{\frac{1}{2}}(\Gamma)$, we obtain \\

\begin{equation*}
\begin{split}
\dfrac{1}{2}\dfrac{d}{dt}|k^{\varepsilon}-g|^{2}_{H^{\frac{1}{2}}(\Gamma)} +& \dfrac{1}{\varepsilon}|k^{\varepsilon}-g|^{2}_{H^{\frac{1}{2}}(\Gamma)}=-(g_{t},k^{\varepsilon}-g)\\
&\leq |g_{t}|_{H^{\frac{1}{2}}(\Gamma)}|k^{\varepsilon}-g|_{H^{\frac{1}{2}}(\Gamma)}\\
&\leq \dfrac{\varepsilon}{2}|g_{t}|^{2}_{H^{\frac{1}{2}}(\Gamma)}+\dfrac{1}{2\varepsilon}|k^{\varepsilon}-g|^{2}_{H^{\frac{1}{2}}(\Gamma)}.
\end{split}
\end{equation*}

\vspace{0.03cm} Hence
\begin{equation}\label{e2.6}
\dfrac{d}{dt}|k^{\varepsilon}-g|^{2}_{H^{\frac{1}{2}}(\Gamma)}+\frac{1}{\varepsilon}|k^{\varepsilon}-g|^{2}_{H^{\frac{1}{2}}(\Gamma)}\leq \varepsilon |g_{t}|^{2}_{H^{\frac{1}{2}}(\Gamma)}.
\end{equation}

\vspace{0.03cm} Using the Gronwall inequality we obtain\\
 \begin{equation}\label{e2.7}
 |k^{\varepsilon}-g|^{2}_{H^{\frac{1}{2}}(\Gamma)}(t)\leq e^{-\frac{t}{\varepsilon}}|k^{\varepsilon}-g|^{2}_{H^{\frac{1}{2}}(\Gamma)}(0)+\varepsilon|g_{t}|^{2}_{L^{2}(0,T;H^{\frac{1}{2}}(\Gamma))}.
 \end{equation}

 \vspace{0.03cm} We integrate (\ref{e2.7}) over $(0,T)$, and we obtain ($u_{0}\in H^{1}(\Omega)\
 \footnote{We do not address the question of minimal regularity of $u_0$, that is
 e.g.~$ u_0\in L^2(\Omega),$ which is not in the scope of this article.
 Indeed the problem of incompatible data occurs already with very smooth data.}$)

\begin{equation}\label{e2.7a}
\int^{T}_{0}|k^{\varepsilon}-g|^{2}_{H^{\frac{1}{2}}(\Gamma)}dt\leq\varepsilon(1-e^{-\frac{T}{\varepsilon}})\Bigl|u_{0}|_{\Gamma}-g(0)\Bigr|^{2}_{H^{\frac{1}{2}}(\Gamma)}
+\varepsilon T|g_{t}|^{2}_{L^{2}(0,T;H^{\frac{1}{2}}(\Gamma))},
\end{equation}
and hence
\begin{equation*}
|k^{\varepsilon}-g|_{L^{2}(0,T;H^{\frac{1}{2}}(\Gamma))}=O(\sqrt{\varepsilon}),
\end{equation*}
which implies
\begin{equation}\label{e2.8}
 k^{\varepsilon}\rightarrow g \enspace strongly\enspace
  in\enspace
 L^{2}(0,T;H^{\frac{1}{2}}(\Gamma))\enspace as\enspace\varepsilon\rightarrow 0.
 \end{equation}

Integrating (\ref{e2.5}) from 0 to t, we obtain
\begin{equation}\label{e2.9}
\int^{t}_{0}(k^{\varepsilon}-g)ds = -\varepsilon (k^{\varepsilon}-g)(t)-\varepsilon g(t) +\varepsilon u_{0}|_{\partial\Omega},
\end{equation}
which yields
\begin{equation}\label{e2.10}
\int^{t}_{0}k^{\varepsilon}(s)ds\rightarrow \int^{t}_{0}g(s)ds\enspace strongly\enspace in\enspace L^{2}(0,T;H^{\frac{1}{2}}(\Gamma))\enspace as\enspace\varepsilon\rightarrow 0.
\end{equation}
The proof of Lemma \ref{l2.1} is complete.

\end{proof}

\begin{remark}\label{r2.2}
We could prove a stronger result namely,
$\displaystyle k^{\varepsilon}\rightarrow g,\int^{t}_{0}k^{\varepsilon}(s)ds$ $\displaystyle \rightarrow \int^{t}_{0}g(s)ds$ strongly in $L^{q}(0,T;H^{\frac{1}{2}}(\Gamma))$, for all $1\leq q < \infty$. But in this article,  $q=2$ is enough for our needs. And for $q=\infty$, from (\ref{e2.7}), we obtain
\begin{equation}\label{e2.10a}
k^{\varepsilon}-g=O(\sqrt{\varepsilon})\enspace in\enspace L^{\infty}(t_{0},T;H^{\frac{1}{2}}(\Gamma))\enspace for\enspace\forall t_{0}>0,
\end{equation}
and also from (\ref{e2.9}), because $k^{\varepsilon}-g$ and $g$ are bounded in $L^{\infty}(0,T;H^{\frac{1}{2}}(\Gamma))$, we obtain
\begin{equation}\label{e2.10b}
\int^{t}_{0}k^{\varepsilon}(s)ds\rightarrow \int^{t}_{0}g(s)ds\enspace strongly\enspace in\enspace L^{\infty}(0,T;H^{\frac{1}{2}}(\Gamma)),
\end{equation}
the norm of the difference being of order $\varepsilon$.
\end{remark}

\subsection{Convergence results for $u^{\varepsilon}$}
Since $\Omega$ is smooth, there exists a lifting operator $L$, linear continuous from $H^{\frac{1}{2}}(\Gamma)$ to $ H^{1}(\Omega)$. We consider such an operator and set $K^{\varepsilon}=L(k^{\varepsilon}), G=L(g)$, and thus have by assumption $G\in L^{\infty}(0,T;H^{1}(\Omega))$, $G_{t}\in L^{2}(0,T;H^{1}(\Omega))$. So we immediately infer from (\ref{e2.8}), (\ref{e2.10}), (\ref{e2.10a}) and (\ref{e2.10b}) that, as $\varepsilon\rightarrow 0$
\begin{equation}\label{e2.11}
 K^{\varepsilon}\rightarrow G \enspace strongly\enspace in\enspace
 L^{2}(0,T;H^{1}(\Omega))\cap L^{\infty}(t_{0},T;H^{1}(\Omega))\enspace\forall t_{0}>0,
 \end{equation}
\begin{equation}\label{e2.12}
\int^{t}_{0}K^{\varepsilon}(s)ds\rightarrow \int^{t}_{0}G(s)ds \enspace strongly\enspace in\enspace L^{\infty}(0,T;H^{1}(\Omega)).
\end{equation}
We now prove Theorem \ref{t1}.\\

\textbf{Proof of Theorem \ref{t1}.} Set $v^{\varepsilon}=u^{\varepsilon}-K^{\varepsilon}$; then
the system (\ref{e2.1}) yields\\
\begin{equation}\left\{\label{e2.17}
\begin{aligned}
&v^{\varepsilon}_{t}-\nu \triangle v^{\varepsilon}=f-K^{\varepsilon}_{t}+\nu\triangle K^{\varepsilon},\\
&v^{\varepsilon}|_{t=0}=u_{0}-K^{\varepsilon}(0)=u_{0}-Lu_{0}|_{\p\Omega},\\
&v^{\varepsilon}|_{\partial\Omega}=0.
\end{aligned}\right.
\end{equation}
Integrating $(\ref{e2.17})_{1}$ from 0 to $t$, we obtain
\begin{equation}\label{e2.17a}
v^{\varepsilon}(t)-\nu\triangle\int^{t}_{0}v^{\varepsilon}(s)ds=\int^{t}_{0}f(s)ds-K^{\varepsilon}(t)+\nu\triangle\int^{t}_{0}K^{\varepsilon}(s)ds+u_{0}.
\end{equation}

We set $\displaystyle V^{\varepsilon}=\int^{t}_{0}v^{\varepsilon}(s)ds$ (with$V^{\varepsilon}(0)=0$), and $\displaystyle F(t)=\int^{t}_{0}f(s)ds+u_{0}$; so $V^{\varepsilon}$ solves the following system
\begin{equation}\left\{\label{e2.17b}
\begin{aligned}
&V^{\varepsilon}_{t}-\nu\triangle V^{\varepsilon}=F-K^{\varepsilon}+\nu\triangle\int^{t}_{0}K^{\varepsilon}(s)ds,\\
&V^{\varepsilon}|_{t=0}=0,\\
&V^{\varepsilon}|_{\partial\Omega}=0.
\end{aligned}\right.
\end{equation}

We take the scalar product of $(\ref{e2.17b})_{1}$ with $V^{\varepsilon}$ in $L^{2}(\Omega)$ and find,
\begin{equation}\label{e2.18}
\frac{1}{2}\frac{d}{dt}|V^{\varepsilon}|^{2}+\nu\parallel V^{\varepsilon}\parallel^{2}=(F,V^{\varepsilon})- (K^{\varepsilon},V^{\varepsilon})+\nu(\triangle\int^{t}_{0} K^{\varepsilon}(s)ds,V^{\varepsilon}).
\end{equation}

We can bound the terms in the right-hand-side of (\ref{e2.18}) as follows:
\begin{equation}\label{e2.19}
(F,V^{\varepsilon})\leq |F||V^{\varepsilon}|\leq c_{1}|F|\parallel V^{\varepsilon}\parallel\leq c_{1}'|F|^{2}+\frac{\nu}{6}\parallel V^{\varepsilon}\parallel^{2},
\end{equation}
\begin{equation}\label{e2.20}
-(K^{\varepsilon},V^{\varepsilon}) \leq
|K^{\varepsilon}||V^{\varepsilon}|\leq
c_{1}|K^{\varepsilon}|\parallel V^{\varepsilon}\parallel\leq
c_{2}'|K^{\varepsilon}|^{2}+\dfrac{\nu}{6}\parallel
V^{\varepsilon}\parallel^{2},
\end{equation}

\begin{equation}\label{e2.21}
\begin{aligned}
\nu(\triangle\int^{t}_{0} K^{\varepsilon}(s)ds,V^{\varepsilon}) &= -\nu(\nabla\int^{t}_{0}K^{\varepsilon}(s)ds,\nabla V^{\varepsilon})\\
&\leq\nu\parallel \int^{t}_{0}K^{\varepsilon}(s)ds\parallel\parallel V^{\varepsilon}\parallel\\
&\leq\dfrac{\nu}{6}\parallel
V^{\varepsilon}\parallel^{2}+c_{3}'\parallel\int^{t}_{0}K^{\varepsilon}(s)ds\parallel^{2}.
\end{aligned}
\end{equation}

Here and below, the $c,c',c_{i},c'_{i}$ are various constants independent of $\varepsilon$, which may be different at different places.\\

\vspace{0.03cm} Combining (\ref{e2.18}), (\ref{e2.19}), (\ref{e2.20}) and (\ref{e2.21}) gives
\begin{equation}\label{e2.22}
\dfrac{d}{dt}|V^{\varepsilon}|^{2} +\nu\parallel
V^{\varepsilon}\parallel^{2}\leq
c_{1}'|F|^{2}+c_{2}'|K^{\varepsilon}|^{2}+c_{3}'\parallel \int^{t}_{0}K^{\varepsilon}(s)ds\parallel^{2}.
\end{equation}

\vspace{0.03cm}
Integrating (\ref{e2.22}) over $(0,t)$, we obtain\\
\begin{equation}\label{e23}
\begin{aligned}
|V^{\varepsilon}(t)|^{2}+\nu\int^{t}_{0}\parallel
V^{\varepsilon}\parallel^{2}ds & \leq c_{1}'\int^{t}_{0}|F|^{2}ds+c_{2}'\int^{t}_{0}|K^{\varepsilon}|^{2}ds\\
&+c_{3}'\int^{t}_{0}\parallel\int^{s}_{0} K^{\varepsilon}(\tau)d\tau\parallel^{2}ds.
\end{aligned}
\end{equation}

We also integrate (\ref{e2.22}) over $(0,T)$, and obtain\\
\begin{equation}\label{e22}
\begin{aligned}
|V^{\varepsilon}(T)|^{2}+\nu\int^{T}_{0}\parallel
V^{\varepsilon}\parallel^{2}ds&\leq c_{1}'\int^{T}_{0}|F|^{2}ds+c_{2}'\int^{T}_{0}|K^{\varepsilon}|^{2}ds\\
&+c_{3}'\int^{T}_{0}\parallel \int^{s}_{0}K^{\varepsilon}(\tau)d\tau\parallel^{2}ds.
\end{aligned}
\end{equation}

It follows from (\ref{e2.11}) and (\ref{e2.12}) that $K^{\varepsilon}$ and $\displaystyle \int^{t}_{0}K^{\varepsilon}(s)ds$ are bounded in $ L^{2}(0,T;H^{1}(\Omega))$, and thus (\ref{e23}), (\ref{e22}) yields:
\begin{equation}\label{e2.25}
 V^{\varepsilon}\enspace remains \enspace bounded\enspace in\enspace L^{\infty}(0,T;L^{2}(\Omega))\cap L^{2}(0,T;H^{1}_{0}(\Omega))\enspace as \enspace \varepsilon\rightarrow 0.
\end{equation}
Thus, there exists a subsequence $V^{\varepsilon'}$ and $V\in L^{\infty}(0,T;L^{2}(\Omega))\cap L^{2}(0,T;$\\
$H^{1}_{0}(\Omega))$ such
that, as $\varepsilon'\rightarrow 0$,
\begin{equation}\label{e2.27}
V^{\varepsilon'}\rightarrow V\enspace weakly\enspace in\enspace L^{2}(0,T;H^{1}_{0}(\Omega)),
\end{equation}

\vspace{0.03cm}\hspace{3cm} $and \enspace weak-star\enspace in\enspace L^{\infty}(0,T;L^{2}(\Omega)).$\\

\vspace{0.03cm} Using (\ref{e2.11}), (\ref{e2.12}) and (\ref{e2.27}), we can pass to the limit in (\ref{e2.17}) with the sequence $\varepsilon'\rightarrow 0$. We proceed as follows.\\

For all $a\in H^{1}_{0}(\Omega)$, and $\phi$ in $C^{1}(0,T)$
 with $\phi(T)=0$,
we multiply $(\ref{e2.17b})_{1}$  by $a\phi$ and integrate over $\Omega\times(0,T)$; we obtain
\begin{equation}\label{e2.29}
\begin{aligned}
-\int^{T}_{0}(V^{\varepsilon'},a)&\phi'(t)dt+\nu\int^{T}_{0}(\nabla V^{\varepsilon'},\nabla a)\phi(t)dt=\int^{T}_{0}(F,a)\phi(t)dt\\
&-\int^{T}_{0}(K^{\varepsilon'},a)\phi(t)dt-\nu\int^{T}_{0}(\nabla \int^{t}_{0}K^{\varepsilon'}(s)ds,\nabla a)\phi(t)dt.
\end{aligned}
\end{equation}

\vspace{0.05cm}Passing to the limit with (\ref{e2.11}), (\ref{e2.12}),
(\ref{e2.27}), we find
\begin{equation}\label{e2.30}
\begin{aligned}
-\int^{T}_{0}(V,a)\phi'(t)dt&+\nu\int^{T}_{0}(\nabla V,\nabla a)\phi(t)dt=\int^{T}_{0}(F,a)\phi(t)dt\\
&-\int^{T}_{0}(G,a)\phi(t)dt-\nu\int^{T}_{0}(\nabla\int^{t}_{0}
G(s)ds,\nabla a)\phi(t)dt.
\end{aligned}
\end{equation}

\vspace{0.05cm}Taking $\phi\in\mathcal{D}(0,T)$, we see that $V$ satisfies
\begin{equation}\label{e2.31}
(V_{t},a)-\nu(\triangle V,a)=(F-G+\nu\triangle \int^{t}_{0}G(s)ds, a), \quad
\forall\enspace a\in H^{1}_{0}(\Omega).
\end{equation}

Now we want to show that $V(0)=0$.\\

\vspace{0.05cm} We classically integrate $(\ref{e2.31})$ times $\phi(t)$ over $(0,T)$ and we obtain:
\begin{equation}\label{e2.32}
\begin{aligned}
-\int^{T}_{0}(V,a)\phi'(t)dt&+\nu\int^{T}_{0}(\nabla V,\nabla a)\phi(t)dt=\int^{T}_{0}(\nu\triangle\int^{t}_{0}G(s)ds,a)\phi(t)dt\\
&+\int^{T}_{0}(F-G,a)\phi(t)dt+(V(0),a)\phi(0).
\end{aligned}
\end{equation}

\vspace{0.03cm}By comparing with (\ref{e2.30}), we find that
\begin{equation}
(V(0),a)\phi(0)=0,
\end{equation}
for every $a\in H^1_0(\Omega)$ and every $\phi\in\mathcal{C}^1([0,T])$ with $\phi(T)=0.$
This implies $V(0)=0$ as
desired. Finally $V$ satisfies
\begin{equation}\left\{\label{e2.33}
\begin{aligned}
&V_{t}-\nu\triangle V=F -G+\nu\triangle \int^{t}_{0}G(s)ds,\\
&V|_{t=0}=0,\\
&V|_{\partial \Omega}=0.
\end{aligned}\right.
\end{equation}

\begin{remark}\label{r2.3}
Furthermore, we could prove that the whole sequence $V^{\varepsilon}\rightarrow V$ weakly in $L^{2}(0,T;H^{1}_{0}(\Omega))$, and weak star in $L^{\infty}(0,T;L^{2}(\Omega))$. Indeed, if not, arguing by contradiction, we could find  a subsequence $\varepsilon_{i}\rightarrow 0$, such that
\begin{equation}\label{e2.34}
\begin{aligned}
V^{\varepsilon_{i}}\nrightarrow V \enspace in\enspace &L^{2}(0,T;H^{1}_{0}(\Omega))\enspace weakly,\\
&L^{\infty}(0,T;L^{2}(\Omega))\enspace weak-star.
\end{aligned}
\end{equation}
Repeating the argument above leading to (\ref{e2.27}), we could extract from $\varepsilon_{i}$ a subsequence $\varepsilon_{i}'$ and find $\bar{V}$ such that, as $\varepsilon_{i}'\rightarrow 0$,

\begin{equation}\label{e2.37}
\begin{aligned}
V^{\varepsilon_{i}'}\rightarrow \bar{V} \enspace in\enspace &L^{2}(0,T;H^{1}_{0}(\Omega))\enspace weakly,\\
&L^{\infty}(0,T;L^{2}(\Omega))\enspace weak-star,
\end{aligned}
\end{equation}
where $\bar{V}$ is the solution of (\ref{e2.33}). But the solution of (\ref{e2.33}) is unique; hence $V=\bar{V}$, and then (\ref{e2.37}) contradicts (\ref{e2.34}).
\end{remark}

Before we finish the proof of the theorem, we now prove the following Lemma.\\

\textbf{Lemma 2.2}\label{l2.2} \textit{Under the assumptions of Theorem \ref{t1}, with $V,V^{\varepsilon}$ being the solutions of (\ref{e2.33}) and (\ref{e2.17b}), we have, as $\varepsilon\rightarrow 0$,}
\begin{equation}\label{e2.35}
V^{\varepsilon}\rightarrow V\enspace strongly \enspace in\enspace L^{2}(0,T;H^{1}_{0}(\Omega))\cap \mathcal{C}([0,T];L^{2}(\Omega)).
\end{equation}

\vspace{0.03cm} \textbf{Proof.} We subtract $(\ref{e2.17b})_{1}$ from $(\ref{e2.33})_{1}$, and obtain
\begin{equation}\label{e2.36}
V_{t}-V^{\varepsilon}_{t}-\nu(\triangle V-\triangle V^{\varepsilon})=K^{\varepsilon}-G+\nu(\triangle\int^{t}_{0} G(s)ds-\triangle \int^{t}_{0}K^{\varepsilon}(s)ds).
\end{equation}

\vspace{0.03cm} We then take the scalar product of (\ref{e2.36}) with $V-V^{\varepsilon}$ in $L^{2}(\Omega)$, and integrate in time from $0$ to $t$, and we obtain:
\begin{equation}\label{e2.36a}
\begin{aligned}
\dfrac{1}{2}|(V-V^{\varepsilon})(t)|^{2} &+\nu\int^{t}_{0}\parallel V-V^{\varepsilon }\parallel^{2}ds=\int^{t}_{0}(K^{\varepsilon}-G,V-V^{\varepsilon})ds\\
&\int^{t}_{0}(\nu\triangle\int^{t}_{0} G(s)ds-\nu\triangle \int^{t}_{0}K^{\varepsilon}(s)ds,V-V^{\varepsilon})ds.
\end{aligned}
\end{equation}

\vspace{0.03cm}Now we set
 \begin{equation}
 \chi_{\varepsilon}(t)=\dfrac{1}{2}|(V-V^{\varepsilon})(t)|^{2}+\dfrac{\nu}{2}\int^{t}_{0}\parallel V-V^{\varepsilon}\parallel^{2}ds,
  \end{equation}
and estimate the right-hand-side of (\ref{e2.36a}) as follows:
\begin{equation}\label{e2.39}
\begin{aligned}
\int^{t}_{0}(K^{\varepsilon}-G,V-V^{\varepsilon})ds&\leq\int^{t}_{0}|K^{\varepsilon}-G||V-V^{\varepsilon}|ds\\
&\leq c_1\int^{t}_{0}|K^{\varepsilon}-G|\parallel V-V^{\varepsilon}\parallel ds\\
&\leq c_1(\int^{T}_{0}|K^{\varepsilon}-G|^{2}ds)^{\frac{1}{2}}(\int^{t}_{0}\parallel V-V^{\varepsilon}\parallel^{2}ds)^{\frac{1}{2}}\\
&\leq c\int^{T}_{0}|K^{\varepsilon}-G|^{2}ds+\dfrac{\nu}{4}\int^{t}_{0}\parallel V-V^{\varepsilon}\parallel^{2}ds,
\end{aligned}
\end{equation}

\begin{equation}\label{e2.40}
\begin{aligned}
\int^{t}_{0}\nu&(\triangle\int^{s}_{0} G(\tau)d\tau-\triangle\int^{s}_{0} K^{\varepsilon}(\tau)d\tau,V-V^{\varepsilon})ds\\
&\leq \nu\int^{t}_{0}\parallel\int^{s}_{0} (G(\tau)-K^{\varepsilon}(\tau))d\tau\parallel\parallel V-V^{\varepsilon}\parallel ds\\
&\leq c'\int^{T}_{0}\parallel \int^{s}_{0}(G- K^{\varepsilon})(\tau)d\tau\parallel^{2}ds+\dfrac{\nu}{4}\int^{t}_{0}\parallel V-V^{\varepsilon}\parallel^{2}ds.
\end{aligned}
\end{equation}

\vspace{0.03cm}Combining (\ref{e2.39}) and (\ref{e2.40}), we see that
\begin{equation}\label{e2.41}
\begin{aligned}
\chi_{\varepsilon}(t)&\leq c'\int^{T}_{0}\parallel\int^{s}_{0} (G- K^{\varepsilon})(\tau)d\tau\parallel^{2}ds+c\int^{T}_{0}|K^{\varepsilon}-G|^{2}ds\\
\end{aligned}
\end{equation}

The right-hand side of (\ref{e2.41}) converges to $0$ as $\epsilon$ converges to $0$, because of (\ref{e2.11}) and (\ref{e2.12}), and so does $\chi_{\varepsilon}(t)$. For $t=T$, we find
\begin{equation}\label{e2.42}
V^{\varepsilon}\rightarrow V\enspace strongly \enspace in\enspace L^{2}(0,T;H^{1}_{0}(\Omega))\enspace as \enspace \varepsilon\rightarrow 0,
\end{equation}
and taking the supreme of (\ref{e2.41}) with respect to $t$, we see that
\begin{equation}\label{e2.49}
V^{\varepsilon}\rightarrow V \enspace strongly \enspace in\enspace L^{\infty}(0,T;L^{2}(\Omega)) \enspace as \enspace\varepsilon\rightarrow 0.
\end{equation}

The lemma is proved.\qed\\

Now we apply Lemma \ref{r2.2} and obtain as $\varepsilon\rightarrow 0$,
\begin{equation}\label{e2.50}
\triangle V^{\varepsilon}\rightarrow\triangle V\enspace strongly\enspace in\enspace L^{2}(0,T;H^{-1}(\Omega))\cap\mathcal{C}([0,T];H^{-2}(\Omega)),
\end{equation}
and from (\ref{e2.12}), we obtain as $\varepsilon\rightarrow 0$,
\begin{equation}\label{e2.51}
\triangle\int^{t}_{0}K^{\varepsilon}(s)ds\rightarrow \triangle\int^{t}_{0}G(s)ds\enspace strongly\enspace in\enspace L^{\infty}(0,T;H^{-1}(\Omega)),
\end{equation}
so after comparing $(\ref{e2.17b})_{1}$ with $(\ref{e2.33})_{1}$, we conclude  that as $\varepsilon\rightarrow 0$,
\begin{equation}\label{e2.52}
\begin{aligned}
V^{\varepsilon}_{t}\rightarrow V_{t}\enspace strongly\enspace in \enspace &L^{2}(0,T;H^{-1}(\Omega))\enspace and\\
& \mathcal{C}([t_{0},T];H^{-2})\enspace for\enspace \forall t_{0}>0.
\end{aligned}
\end{equation}
Now we define $v=V_{t}$, and take the derivative on $(\ref{e2.33})_{1}$, we obtain that $v$ solves the following system,
\begin{equation}\left\{\label{e2.53}
\begin{aligned}
&v_{t}-\nu\triangle v=f-G_{t}+\nu\triangle G,\\
&v|_{t=0}=u_{0}-G(0),\\
&v|_{\partial\Omega}=0.
\end{aligned}\right.
\end{equation}
So (\ref{e2.52}) yields as $\varepsilon\rightarrow 0$,
\begin{equation}\label{e2.54}
\begin{aligned}
v^{\epsilon}\rightarrow v\enspace strongly\enspace in \enspace &L^{2}(0,T;H^{-1}(\Omega))\enspace and\\
& \mathcal{C}([t_{0},T];H^{-2})\enspace for\enspace \forall t_{0}>0.
\end{aligned}
\end{equation}
The final stage of the proof of Theorem \ref{t1} consists in reinterpreting the results above
that is (\ref{e2.11}), (\ref{e2.12}) and (\ref{e2.54}) in terms of the
convergence of $u^{\varepsilon}=v^{\varepsilon}+K^{\varepsilon}$ towards $u=v+G$; we obtain precisely (\ref{e2.3}).
Theorem \ref{t1} is proven.\qed

\begin{remark}\label{r2.4}
Similarly as Remark \ref{r2.2}, we see that we also have $u^{\varepsilon}\rightarrow u$ strongly in $L^{q}(0,T;H^{-1}(\Omega))$ for all $1\leq q <\infty$.
\end{remark}

\subsection{Boundary layer analysis for $k^{\varepsilon}$}\label{s2.2} In the previous section, Lemma \ref{l2.1} stated that under our assumptions, $k^{\varepsilon}$ strongly converges to $g$ in $L^{2}(0,T;H^{\frac{1}{2}}(\Gamma))$, as $\varepsilon\rightarrow 0$.  Here in order to better compare $k^{\varepsilon}$ and $g$, we are going to study the boundary layer for the system (\ref{e2.2}).\\

Along the asymptotic analysis, we define the outer expansion $\displaystyle k^{\varepsilon}\sim\sum^{\infty}_{j=0}\varepsilon^{j}k^{j}$. By formal identification at each power of $\varepsilon$, we obtain
\begin{equation}\label{e2.2.1}
\begin{aligned}
&O(\varepsilon^{-1}): \qquad k^{0}=g,\\
&O(\varepsilon^{j}):\qquad k^{j}_{t}+k^{j+1}=0,\enspace \forall j\geq 0.
\end{aligned}
\end{equation}

By explicit calculations, we find:
\begin{equation}\label{e2.2.2}
k^{j}=(-1)^{j}g^{(j)},\enspace \forall j\geq 0.
\end{equation}

It is clear that the functions $k^{j}$ of the outer expansion do not generally satisfy the initial condition in (\ref{e2.2}) in the case of interest here where $g(0)\neq u_{0}|_{\partial\Omega}$. To account for this discrepancy, we classically introduce the inner expansion $k^{\varepsilon}\sim\sum^{\infty}_{j=0}\varepsilon^{j}\theta^{j}$, where $\theta^{j}=\theta^{j}(\overline{t}), (\overline{t}=t/\varepsilon)$. Then we find\\

\hspace{3cm}$\displaystyle \sum^{\infty}_{j=0}\varepsilon^{j}\frac{d\theta^{j}}{d\overline{t}}+\sum^{\infty}_{j=0}\varepsilon^{j}\theta(\overline{t})=0$.\\

By formal identification at each power of $\varepsilon$, we obtain the following equations:
\begin{equation}\label{e2.2.3}
\frac{d\theta^{j}}{d\overline{t}}+\theta^{j}(\overline{t})=0, \enspace for \enspace j\geq 0.
\end{equation}

The initial conditions that we choose are:
\begin{equation}\label{e2.2.4}
\begin{aligned}
&\theta^{0}(0)=u_{0}|_{\partial\Omega}-g(0),\\
&\theta^{j}(0)=-k^{j}(0)=(-1)^{j+1}g^{(j)}(0),\enspace for\enspace j\geq 1.
\end{aligned}
\end{equation}

By explicit calculations, we obtain:
\begin{equation}\label{e2.2.5}
\theta^{j}=e^{-\frac{t}{\varepsilon}}\theta^{j}(0),\enspace for \enspace j\geq 0.
\end{equation}

To obtain the asymptotic error estimate, we set
\begin{equation}\label{e2.2.6}
w_{\varepsilon n}=k^{\varepsilon}-k_{\varepsilon n}-\theta_{\varepsilon n},
\end{equation}
where\\

\hspace{3cm}$\displaystyle k_{\varepsilon n}=\sum_{j=0}^{n}\varepsilon^{j}k^{j},\enspace \theta_{\varepsilon n}=\sum^{n}_{j=0}\varepsilon^{j}\theta^{j}$.\\

Now we can conclude as follows.\\

\textbf{Theorem 2.2} \label{t2} \textit{If $g^{(n+1)}\in L^{2}(0,T;H^{\frac{1}{2}}(\Gamma))$ for $n\geq 0$, and $w_{\varepsilon n}$ is defined in (\ref{e2.2.6}), then as $\varepsilon\rightarrow 0$}
\begin{equation}\label{e2.2.7}
\begin{aligned}
&w_{\varepsilon n}=O(\varepsilon^{n+1})\enspace in \enspace L^{2}(0,T;H^{\frac{1}{2}}(\Gamma)),\\
&w_{\varepsilon n}=O(\varepsilon^{n+\frac{1}{2}})\enspace in \enspace L^{\infty}(0,T;H^{\frac{1}{2}}(\Gamma)).
\end{aligned}
\end{equation}

\textbf{Proof}. We firstly notice that $w_{\varepsilon n}$ vanishes at $t=0$. We insert then (\ref{e2.2.6}) into (\ref{e2.2}), and we find:
\begin{equation}\left\{\label{e2.2.8}
\begin{aligned}
&\varepsilon (w_{\varepsilon n})_{t}+w_{\varepsilon n}=\varepsilon^{n+1}(-1)^{n+1}g^{(n+1)},\\
&w_{\varepsilon n}|_{t=0}=0.
\end{aligned}\right.
\end{equation}
We take the $H^{\frac{1}{2}}(\Gamma)$ scalar product of $(\ref{e2.2.8})_{1}$ with $w_{\varepsilon n}$ and integrate over $[0,t]$; we obtain
\begin{equation*}
\begin{aligned}
\frac{\varepsilon}{2}|w_{\varepsilon n}(t)|^{2}_{H^{\frac{1}{2}}(\Gamma)}&+\int^{t}_{0}|w_{\varepsilon n}(s)|^{2}_{H^{\frac{1}{2}}(\Gamma)}ds=\int^{t}_{0}((-\varepsilon)^{n+1}g^{(n+1)},w_{\varepsilon n})_{H^{\frac{1}{2}}(\Gamma)}ds\\
&\leq\frac{1}{2}\int^{t}_{0}|w_{\varepsilon n}(s)|^{2}_{H^{\frac{1}{2}}(\Gamma)}ds+\frac{\varepsilon^{2(n+1)}}{2}\int^{t}_{0}|g^{(n+1)}(s)|^{2}_{H^{\frac{1}{2}}(\Gamma)}ds,
\end{aligned}
\end{equation*}
\begin{equation}\label{e2.2.9}
\varepsilon|w_{\varepsilon n}(t)|^{2}_{H^{\frac{1}{2}}(\Gamma)}+\int^{t}_{0}|w_{\varepsilon n}(s)|^{2}_{H^{\frac{1}{2}}(\Gamma)}ds\leq \varepsilon^{2(n+1)}\int^{T}_{0}|g^{(n+1)}(s)|^{2}_{H^{\frac{1}{2}}(\Gamma)}ds.
\end{equation}

If we set $t=T$ in (\ref{e2.2.9}), we obtain $w_{\varepsilon n}=O(\varepsilon^{n+1})$ in $L^{2}(0,T;H^{\frac{1}{2}}(\Gamma))$, and if we take the supremum of (\ref{e2.2.9}) over $[0,T]$, we obtain $w_{\varepsilon n}=O(\varepsilon^{n+\frac{1}{2}})$ in $L^{\infty}(0,T;H^{\frac{1}{2}}(\Gamma))$.

Theorem \ref{t2} has been proved.\qed
\begin{remark}\label{r2.2a}
If we additionally assume that $g_{tt}\in L^{2}(0,T;H^{\frac{1}{2}})$ in Theorem \ref{t1}, from (\ref{e2.2.7}), setting $n=1$, we find \\

\hspace{0.5cm}$w_{\varepsilon 1}=k^{\epsilon}-g+\epsilon g_{t}-(u_{0}|_{\partial\Omega}-g(0))e^{-t/\varepsilon}+\varepsilon g_{t}(0)e^{-t/\varepsilon}=O(\varepsilon^{3/2})$,\\
in $L^{\infty}(0,T;H^{\frac{1}{2}}(\Gamma))$. Then for any $t_{0}>0$, \\

\hspace{0.5cm}$(k^{\varepsilon}-g)(t_{0})=-\varepsilon g_{t}(t_{0})+O(\varepsilon^{3/2})+e.s.t.$,\\
where $e.s.t.$ means exponentially small term (for all $H^{m}$-norms).

Hence $|(k^{\varepsilon}-g)(t_{0})|_{H^{\frac{1}{2}}}=O(\varepsilon)$, for $\forall t_{0}>0$ fixed. We now do similarly as (\ref{e2.5})-(\ref{e2.7}) for $k^{\varepsilon}_{t}-g_{t}$ and integrate from $t_{0}$ to $t$,
\begin{equation}\label{e2.10c}
 |k^{\varepsilon}_{t}-g_{t}|^{2}_{H^{\frac{1}{2}}(\Gamma)}(t)\leq e^{-\frac{t-t_{0}}{\varepsilon}}|k^{\varepsilon}_{t}-g_{t}|^{2}_{H^{\frac{1}{2}}(\Gamma)}(t_{0})+\varepsilon|g_{tt}|^{2}_{L^{2}(0,T;H^{\frac{1}{2}}(\Gamma))},
 \end{equation}
which yields
\begin{equation}\label{e2.10d}
k^{\varepsilon}_{t}-g_{t}=O(\sqrt{\varepsilon})\enspace in\enspace L^{\infty}(t_{0},T;H^{\frac{1}{2}}(\Gamma)),\enspace \forall t_{0}>0.
\end{equation}

\end{remark}

\section{Numerical Results For the penalty method}\label{s3}
\subsection{Approximations of $k^\varepsilon$\label{s3.1a}} In order to test the efficiency of the proposed penalty method, we will provide in this section and in the next one (Sec. \ref{s3.2})
some numerical results for system (\ref{e2.2}) with $\Omega=(0,1)\times (0,1)$, and $0\leq t\leq 1$; we set $g(t)=sin(t)$ for all $(x,y)\in \partial \Omega$, $u_0(x,y)=sin(\frac{5\pi}{4}x+\frac{3\pi}{4})sin(\frac{5\pi}{4}y+\frac{3\pi}{4})$, in which case, we face the discrepancies all along the lines $x=0$ and $y=0$
(but no discrepancy along the parts $x=1$ or $y=1$ of the boundary).\\

We start by testing the quality of the approximation of $k^\varepsilon$ inferred by the boundary layer analysis of Section \ref{s2.2}; that is $\displaystyle k^\varepsilon\sim\sum^n_{j=0}\varepsilon^j
(k^j+\theta^{j}),$ for suitable $\varepsilon's$ and $n's$.

\begin{figure}[ht]
\scalebox{0.8}{\includegraphics{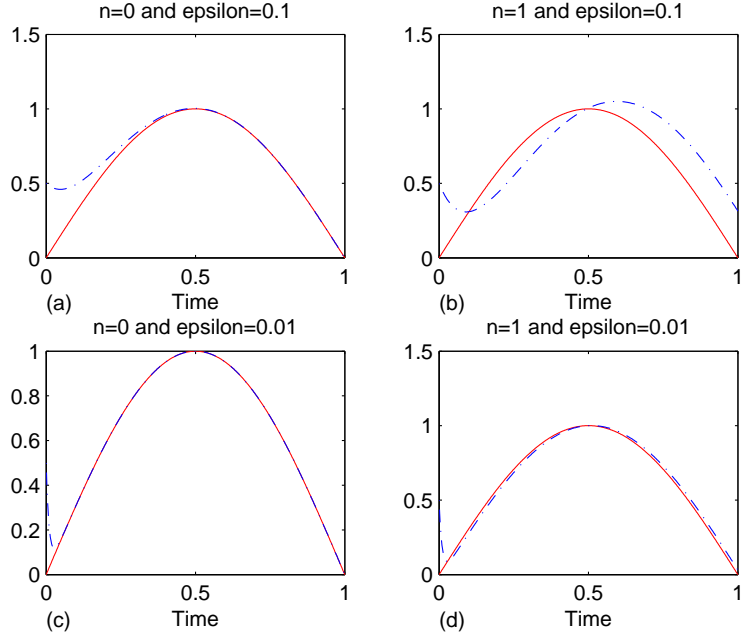}}
\caption{(a) Boundary Layer Element with $\epsilon=0.1,k^{\epsilon}\sim k^{0}+\theta^{0}$, (b) Boundary Layer Element with $\varepsilon=0.1,k^{\varepsilon}\sim k^{0}+\theta^{0}+\varepsilon k^{1}+\varepsilon\theta^{1}$, (c) Boundary Layer Element with $\varepsilon=0.01,k^{\varepsilon}\sim k^{0}+\theta^{0}$, (d) Boundary Layer Element with $\varepsilon=0.01,k^{\varepsilon}\sim k^{0}+\theta^{0}+\varepsilon k^{1}+\varepsilon\theta^{1}$.}
\label{f9}
\end{figure}

Because (\ref{e2.2}) is a $2D$ system which makes the graphing impossible
along the time axis, we restrict ourselves to follow the time evolution of
the exact and approximate function at one point of  the boundary;
for simplicity we choose the corner $(x,y)=(0,0)$. We then plot
in Figure \ref{f9}, $g(t)$ (solid line) and
$\displaystyle k^{\varepsilon}$ (dash--dot line),
$\displaystyle k^\varepsilon\sim\sum^{n}_{j=0}\varepsilon^{j}(k^{j}+\theta^{j})$ with $n=0$ or 1, and $\varepsilon=0.1$ or 0.01.  For $n=1$, the proposed new scheme gives a good approximation of $k^\varepsilon$.\\

\begin{figure}[ht]
\scalebox{0.8}{\includegraphics{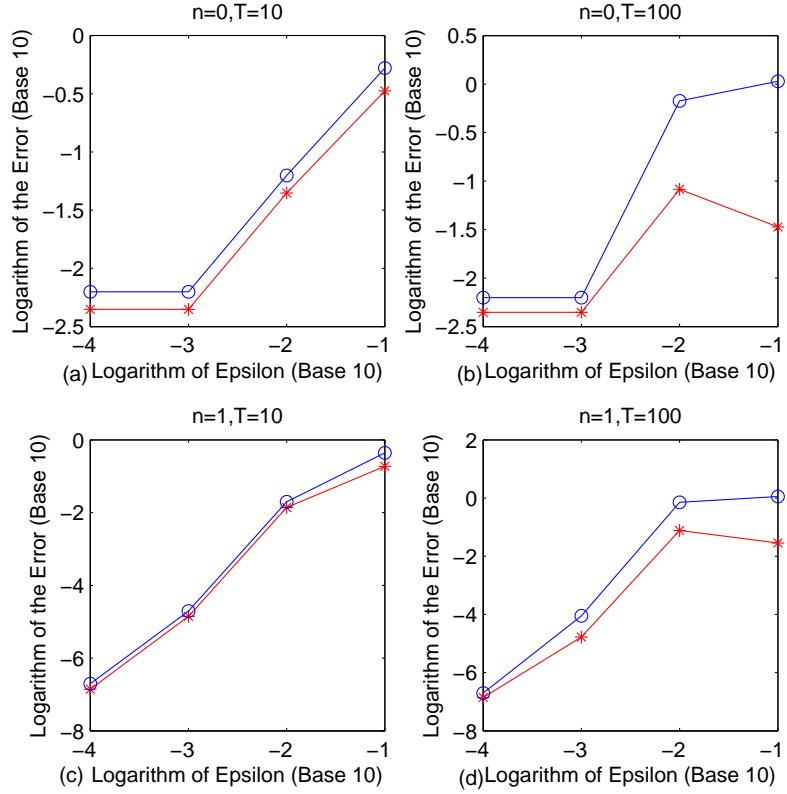}}
\caption{$L^{2}$(plot o)- and $L^{\infty}$(plot *)- error between the boundary layer schemes and the real solution for $\varepsilon k_{t}+k=g$.}
\label{f0.5}
\end{figure}

Fig. \ref{f0.5} gives the $L^{2}-$ and $L^{\infty}-$ errors which stand respectively for the $L^{2}(0,T;H^{\frac{1}{2}}(\Gamma))$ and $L^{\infty}(0,T;H^{\frac{1}{2}}(\Gamma))$
 norms of the difference between the real solution $k^{\varepsilon}$ and the approximations $\displaystyle \sum^{n}_{j=0}\varepsilon^{j}(k^{j}+\theta^{j})$, as the number of time steps $T$, $\varepsilon$ and $n$ vary. It is clear that the smaller $\varepsilon$ is, the smaller both errors are.\\

\begin{figure}[ht]
\scalebox{0.7}{\includegraphics{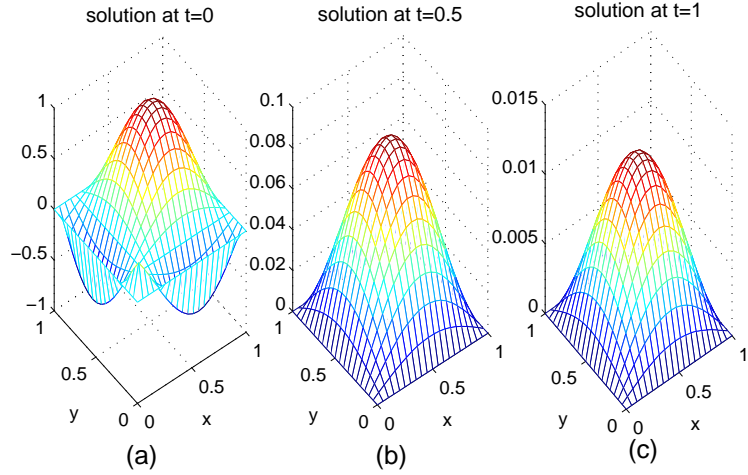}}
\caption{The exact solution of the system in the square $\Omega$ without applying the penalty method, at times $0, 0.5$ and $1$ (Figures \ref{f1}(a), \ref{f1}(b), \ref{f1}(c)).}
\label{f1}
\end{figure}
\begin{figure}[ht]
\scalebox{0.55}{\includegraphics{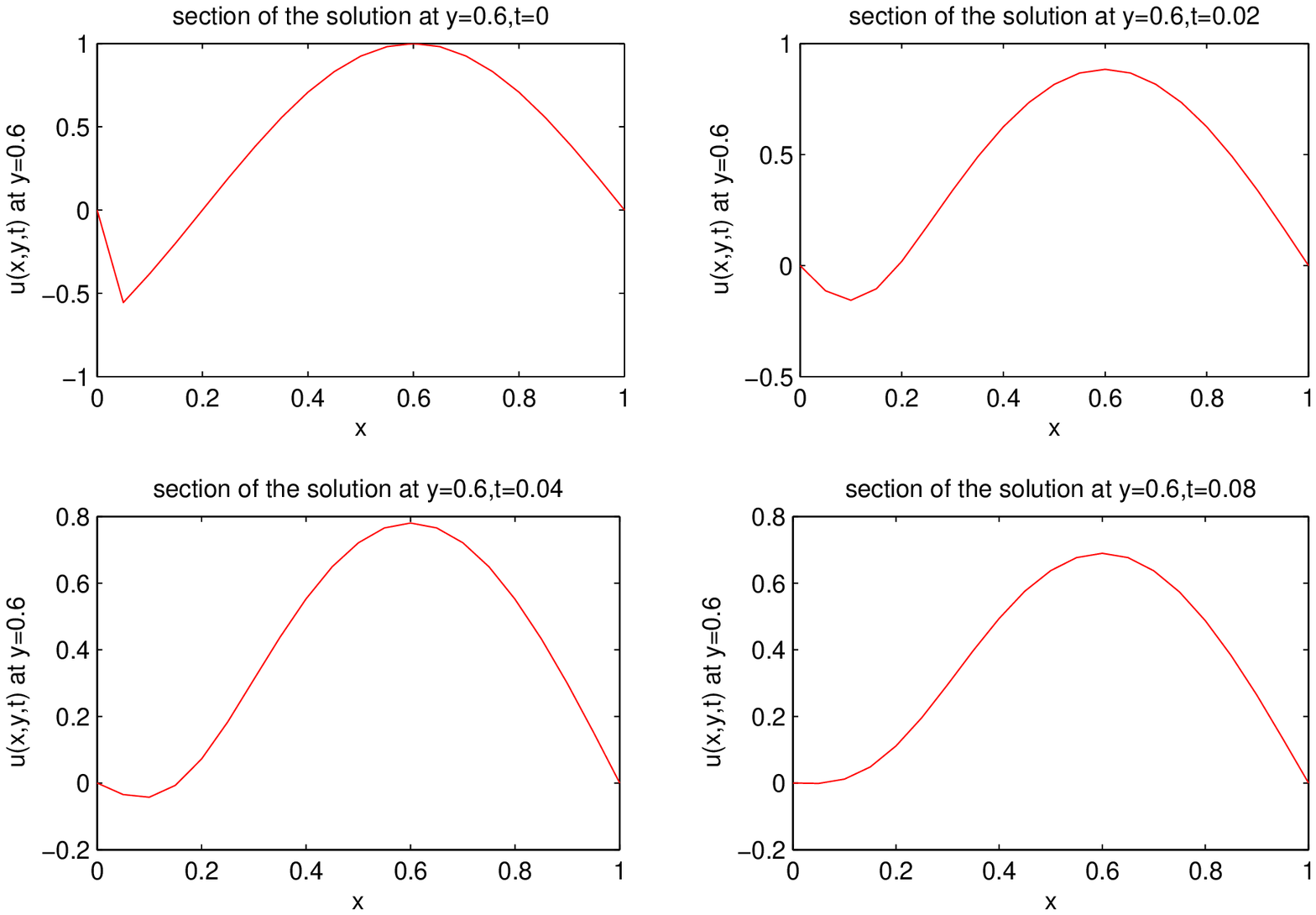}}
\caption{The sections of the exact solution at $y=0.6$ when $t$ is close to $0$.}\label{f8}
\end{figure}

\subsection{A two-dimensional system in a square $\Omega$\label{s3.1}}

To verify the effectiveness of the penalty method, we use the finite
elements method for the spatial approximation of $u.$  The Penalty
Method is mainly aimed for multi-dimensional time-dependent PDEs, so we consider
the $2D$ system, as in (\ref{e1.2}):
\begin{figure}[ht]
\scalebox{0.75}{\includegraphics{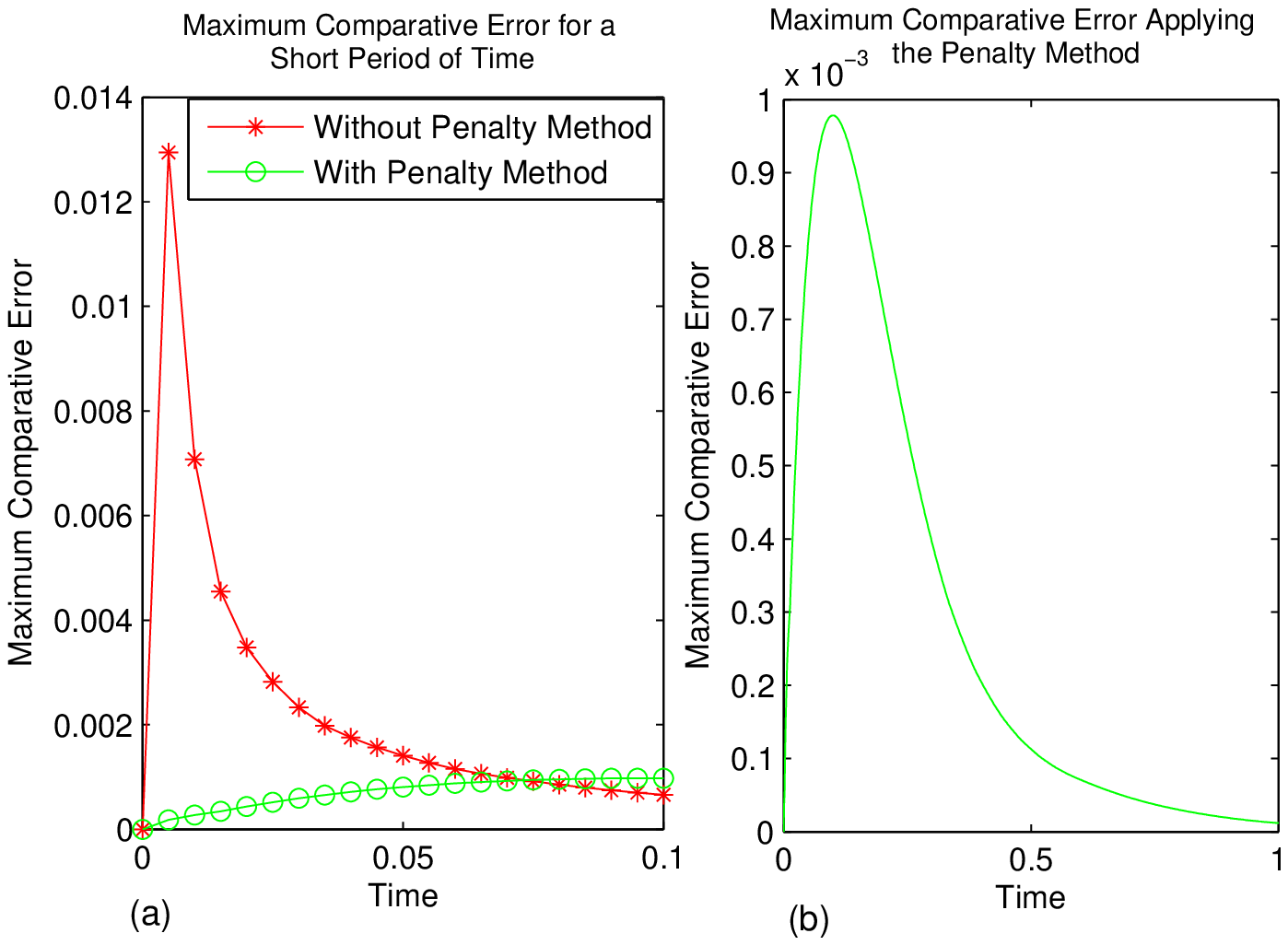}} \caption{The comparative
errors of the 2D system in the $L^{\infty}$ norm for $\varepsilon=0.1$. (a) Maximum comparative error for a short time period (in real value), (b) Maximum comparative error when applying the Penalty
Method (\textit{times} $10^{-3}$).}\label{f2}
\end{figure}
\begin{equation}\left\{\label{e3.1}
\begin{aligned}
&\frac{\partial u}{\partial t}-\nu (u_{xx}+u_{yy})=f,\\
&u|_{\partial \Omega}= g,\\
&u|_{t=0}=u_{0}.
\end{aligned}\right.
\end{equation}
where $0\leq x\leq 1, 0\leq y\leq 1, 0\leq t\leq 1$.\\
\begin{figure}[ht]
\scalebox{0.69}{\includegraphics{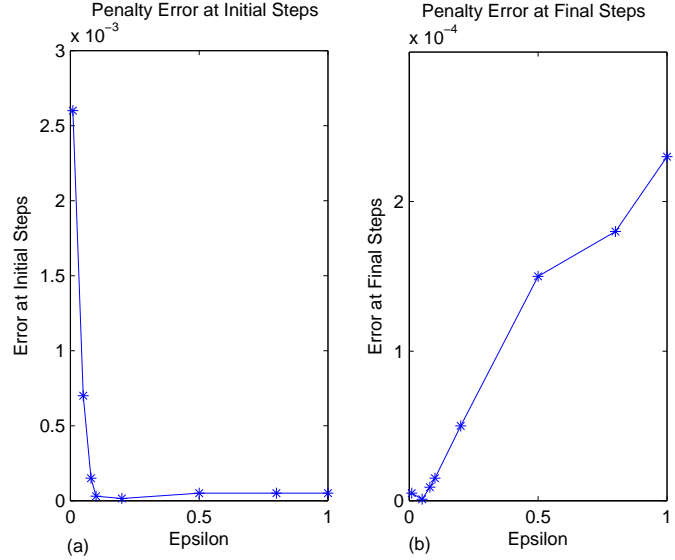}}\caption{The comparative errors for the 2D system in the $L^{\infty}$ norm with variations of $\epsilon$ with mesh $\triangle x=\frac{1}{24},\triangle y=\frac{1}{24},\triangle t=\frac{1}{1000}$. (a) The error at initial step. (b) The error at final step.  Note the factors $10^{-3}, 10^{-4}$ in (a), (b)}\label{f5}
\end{figure}

\begin{figure}[ht]
\scalebox{0.7}{\includegraphics{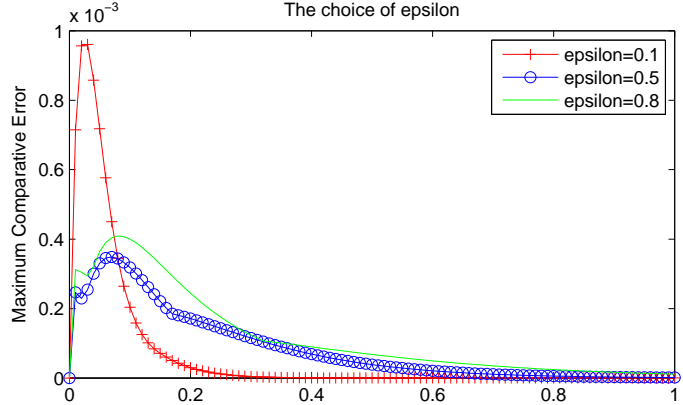}}\caption{The maximum comparative errors for the 2D system in square domain.}\label{f4.1}
\end{figure}

\begin{figure}[ht]
\scalebox{0.56}{\includegraphics{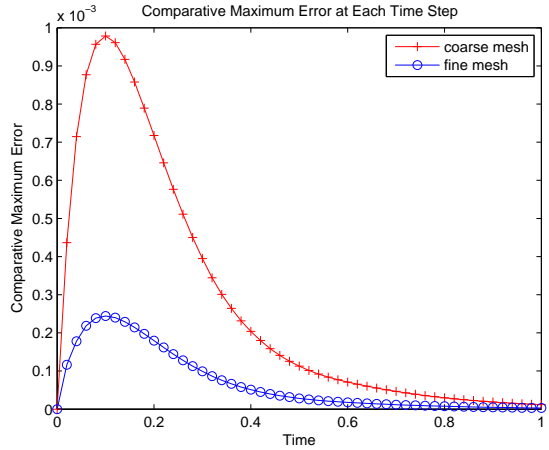}}\caption{The comparative errors for the 2D system in $L^{\infty}$ sense at $\epsilon=0.1$. The upper line is with mesh $\triangle x=\frac{1}{24},\triangle y=\frac{1}{24},\triangle t=\frac{1}{1000}$, the lower line is  with mesh $\triangle x=\frac{1}{48},\triangle y=\frac{1}{48},\triangle t=\frac{1}{4000}$ }\label{f4}
\end{figure}

\begin{figure}[ht]
\scalebox{0.7}{\includegraphics{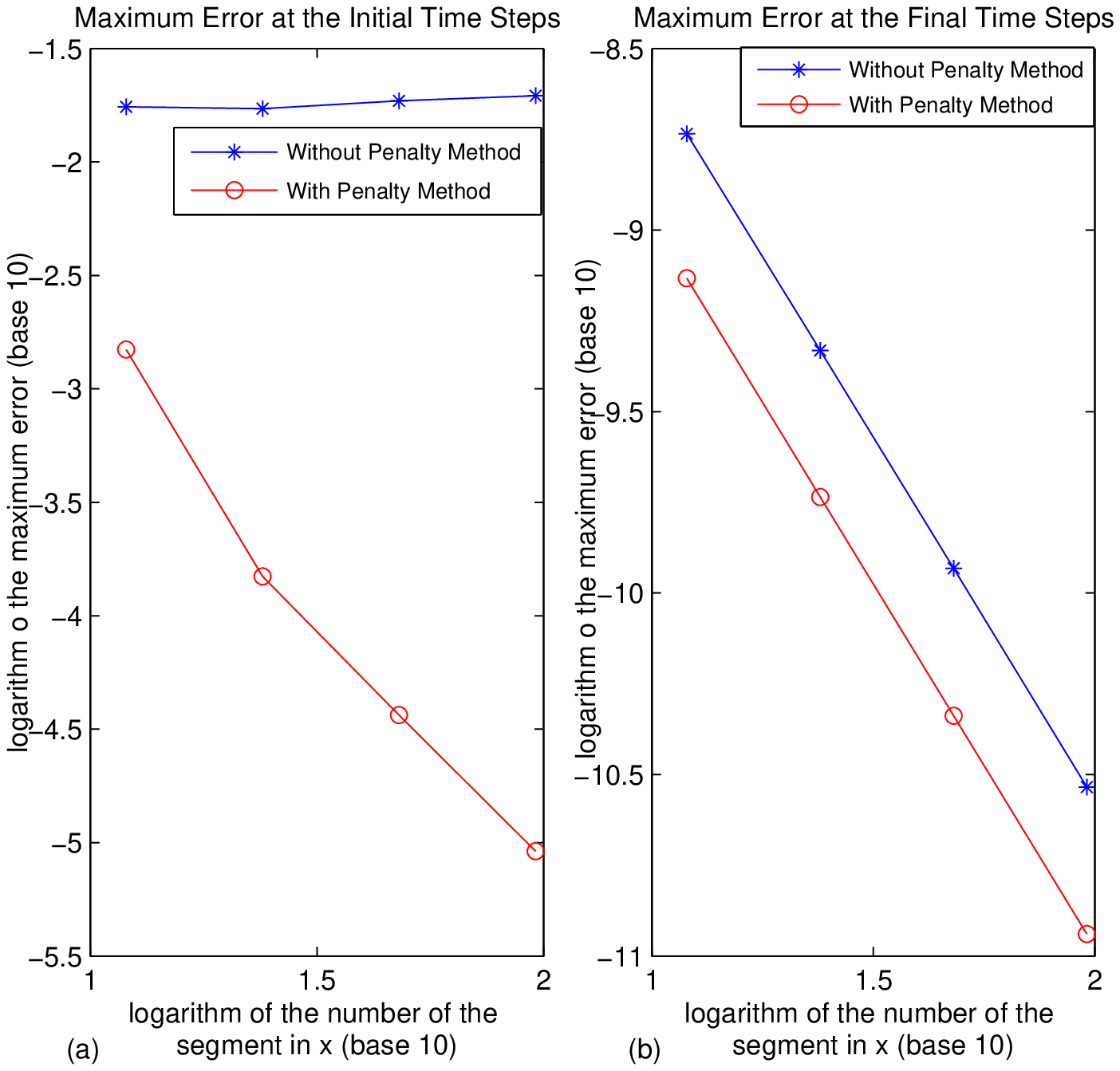}}\caption{Decay of the maximum errors. When we apply the penalty method here $\ep=0.1$. (a) at the initial steps, (b) at the final steps}\label{f4.5}
\end{figure}


\vspace{0.03cm}We set $\nu=0.2$, $f=0,\enspace g=0,\enspace
u_{0}=sin(\dfrac{5\pi}{4}x+\dfrac{3\pi}{4})sin(\dfrac{5\pi}{4}y+\dfrac{3\pi}{4})$,
so that $g(0)\neq u_{0}|_{\partial\Omega}$ on the lines $x=0$ and $y=0$. So we face the
incompatibility problem, namely the boundary conditions and the initial
condition do not match at these corners of the time and
spatial axes. For the test we set $\epsilon =0.1$ in the
penalty approximation \eqref{e2.1}--\eqref{e2.2} of \eqref{e3.1}. In more general problems, we might have discontinuities at the space corners $x=0$ or $1$, $y=0$ or $1$. But since the function $u_{0}$ is smooth at the corners, these singularities do not occur here, at least at the low orders. \\

We first plot the solution of system (\ref{e3.1}) without
applying the penalty method. The solution is plotted in
Fig. \ref{f1}; (a) is the graph of the approximate solution at $t=0$, (b) is
the graph of the approximate solution at $t=0.5$, and (c) is the graph of the
approximate solution at $t=1$. The graph displays a sharp gradient around the corner of
the time--space axis during an initial short period due to the
incompatibility between
the initial and boundary conditions there. In order to see  the sharp
gradient clearly and the changes of the gradient as time evolves,
we plot the sections $(x\in(0,1), y=0.6)$ of the solution at times close to $0$;
see Figure \ref{f8}. It is clear that, at $t=0$, we observe
the sharpest gradient at the time--space corner and as time evolves,
the gradient becomes smoother and smoother at that corner,
until $t=0.08$, when it is essentially flat. \\

Next, to study the accuracy of the numerical method, we must
measure the errors for the approximate solutions.
Hence we compute the comparative
errors which are the differences between two numerical solutions for
the problem, one with the stated mesh sizes, and the other one with
a finer mesh. Then at each time step, we obtain the maximum error between the two meshes above; it is understood to be $L^{\infty}$ comparative errors, or maximum comparative errors. In what follows, all the error terms are to be understood in this sense.\\

We plot the maximum comparative errors of the $2D$ system on Fig.~\ref{f2}.  Graph (b) is
the plot of the maximum comparative errors along the whole time period if
we apply the penalty method. Because the discrepancy happens at the
time-space corner, we zoom into the left corner of graph (b) and compare
it with the error when we do no apply the penalty method. In graph (a),
the line with stars is the maximum comparative errors with the penalty
method applied, and the line with circles is the maximum comparative
errors without the penalty method. We observe that
the magnitude of the errors at the time-space corner is reduced by around one order
by the penalty method.\\

Because we use finite differences method, so for the same $\varepsilon$,
if we have a finer mesh, the maximum comparative error should be smaller.
In Fig. \ref{f4} we plot the error of the $2D$ system with $\varepsilon=0.1$,
the lower curve with a finer mesh, the upper curve with a coarser mesh.
The magnitude of the errors are reduced by around 40\%. So for a
fixed $\varepsilon$, the finer the mesh is, the smaller the error is. We are also interested in the decay of the maximum errors. The most interesting and informative comparison can be made between the decay rates of the maximum errors at the initial and final time steps. In Fig. \ref{f4.5} (a), the maximum errors at the initial time step are plotted against the grid resolution
in the log–log scale. Without the penalty method, the maximum errors do not decrease as the grid refines, which demonstrates that the
singularity in the solution during the initial period is serious. With
the penalty method $(\ep = 0.1)$, the maximum errors decay at roughly the second order. Fig. \ref{f4.5} (b) shows that, with and without applying penalty method, the maximum errors at the final steps (t=1) decay at approximately the second order.\\
\begin{figure}[ht]
\scalebox{0.7}{\includegraphics{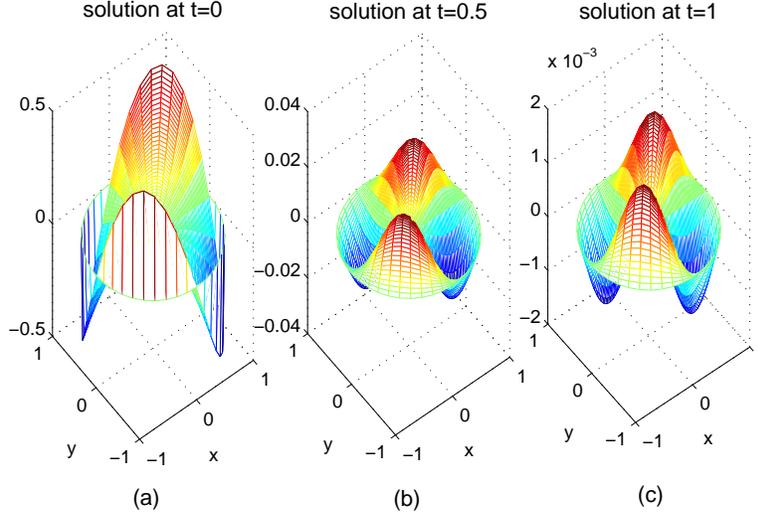}}
\caption{The solution of the system in disk $\Omega$ without applying the penalty method.}
\label{f10}
\end{figure}

Next, we fix the meshes at e.g.
$\triangle x=\frac{1}{24},\triangle y=\frac{1}{24},\triangle t=\frac{1}{1000}$,
and let $\varepsilon$ vary.  In Fig. \ref{f5}, we plot the maximum comparative errors of system (\ref{e3.1}). At the initial steps, the error decreases sharply as $\varepsilon$ increases and remains close to $0$, then it becomes stable flat. At the final steps, the error increases almost linearly as $\varepsilon$ increases. With $\varepsilon$ at about $0.1$, the initial error is minimized while the error at final step is well controlled.
But as Fig. \ref{f4.1} shows, in a short time period $\varepsilon=0.5$ gives us smaller errors and again after a short time period $\varepsilon=0.1$ gives us a smaller error. In optimization theory, the choice of $\varepsilon$ is usually made by trial and error and is not a major issue. It does not follow the "intuitive" idea that the error becomes smaller as $\varepsilon$ becomes smaller because of many other contingent errors such as round-off and descretization errors. In general the choice of $\varepsilon$ depends on our goals of the computation.\\

\subsection{2D system in a disk $\Omega$}\label{s3.2}
To further verify the effectiveness of the penalty method we now test
the results in a different domain.  We now choose a disk $\Omega=\{(x,y)|x^{2}+y^{2}\leq 1\}$.
The 2D heat equations in the polar coordinates
$x=rcos(\theta),y=rsin(\theta)$ where $0\leq \theta\leq 2\pi,0\leq r\leq 1$ read
\begin{equation}\left\{\label{e3.5}
\begin{aligned}
&\frac{\partial u}{\partial t}-\nu(u_{rr}+\frac{u_{r}}{r}+\frac{u_{\theta\theta}}{r^{2}})=f,\\
&u|_{t=0}=u_{0},\\
&u|_{r=1}=g.
\end{aligned}\right.
\end{equation}

\begin{figure}[ht]
\scalebox{0.8}{\includegraphics{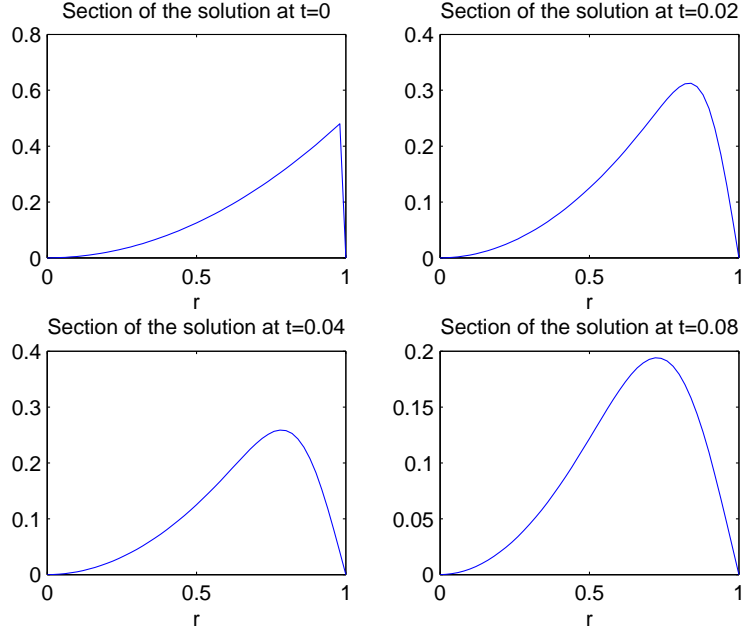}}
\caption{The section of the solution at $\theta=\frac{\pi}{4}$ when t is close to 0.}
\label{f13}
\end{figure}

\begin{figure}[ht]
\scalebox{0.8}{\includegraphics{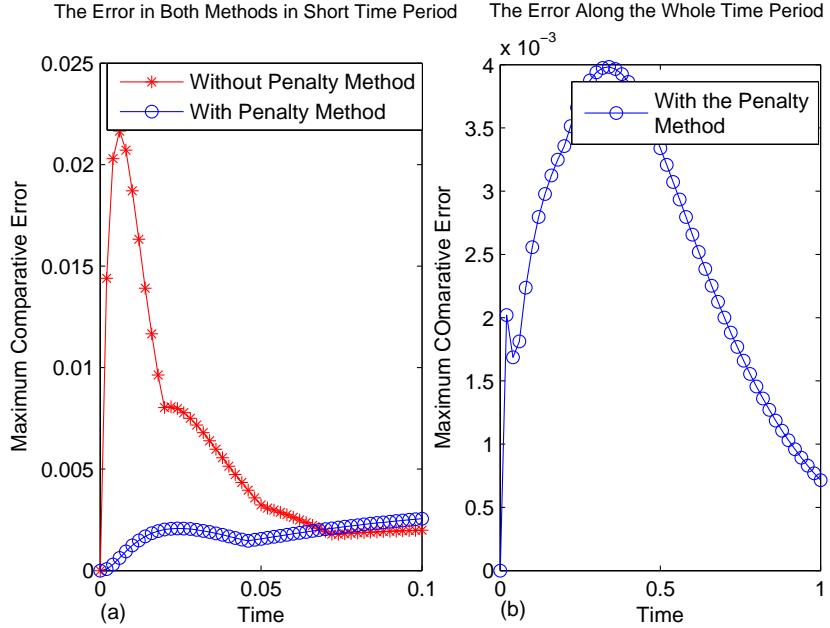}}
\caption{The maximum comparative error for Ordinary Finite Element and penalty method in $L^{\infty}$ sense with $\varepsilon=0.1$}
\label{f11}
\end{figure}

\begin{figure}[ht]
\scalebox{0.8}{\includegraphics{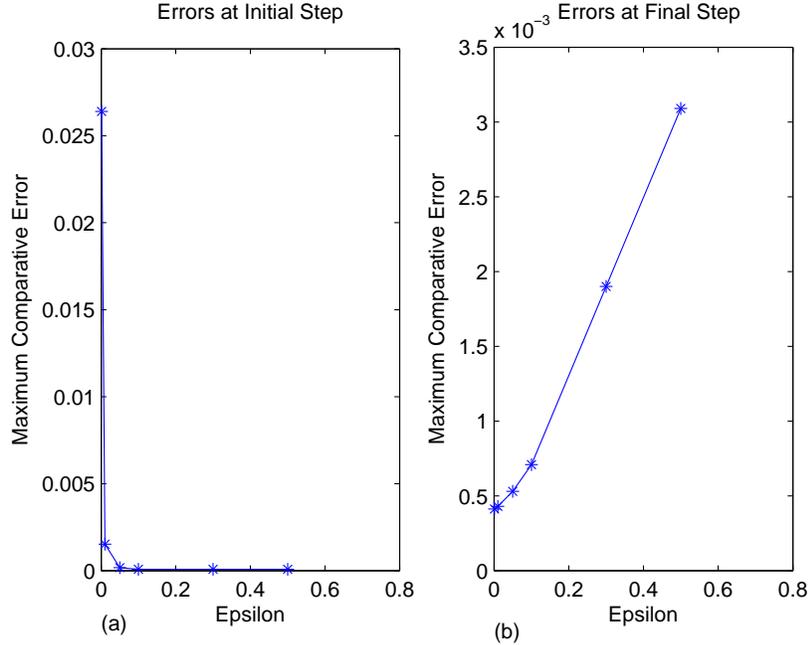}}
\caption{The maximum comparative errors for penalty method at both initial and final steps as $\varepsilon$ variants with mesh $\triangle r=\frac{1}{10},\triangle \theta=\frac{1}{10},\triangle T=\frac{1}{1000}$. }
\label{f12}
\end{figure}
\begin{figure}[ht]
\scalebox{0.7}{\includegraphics{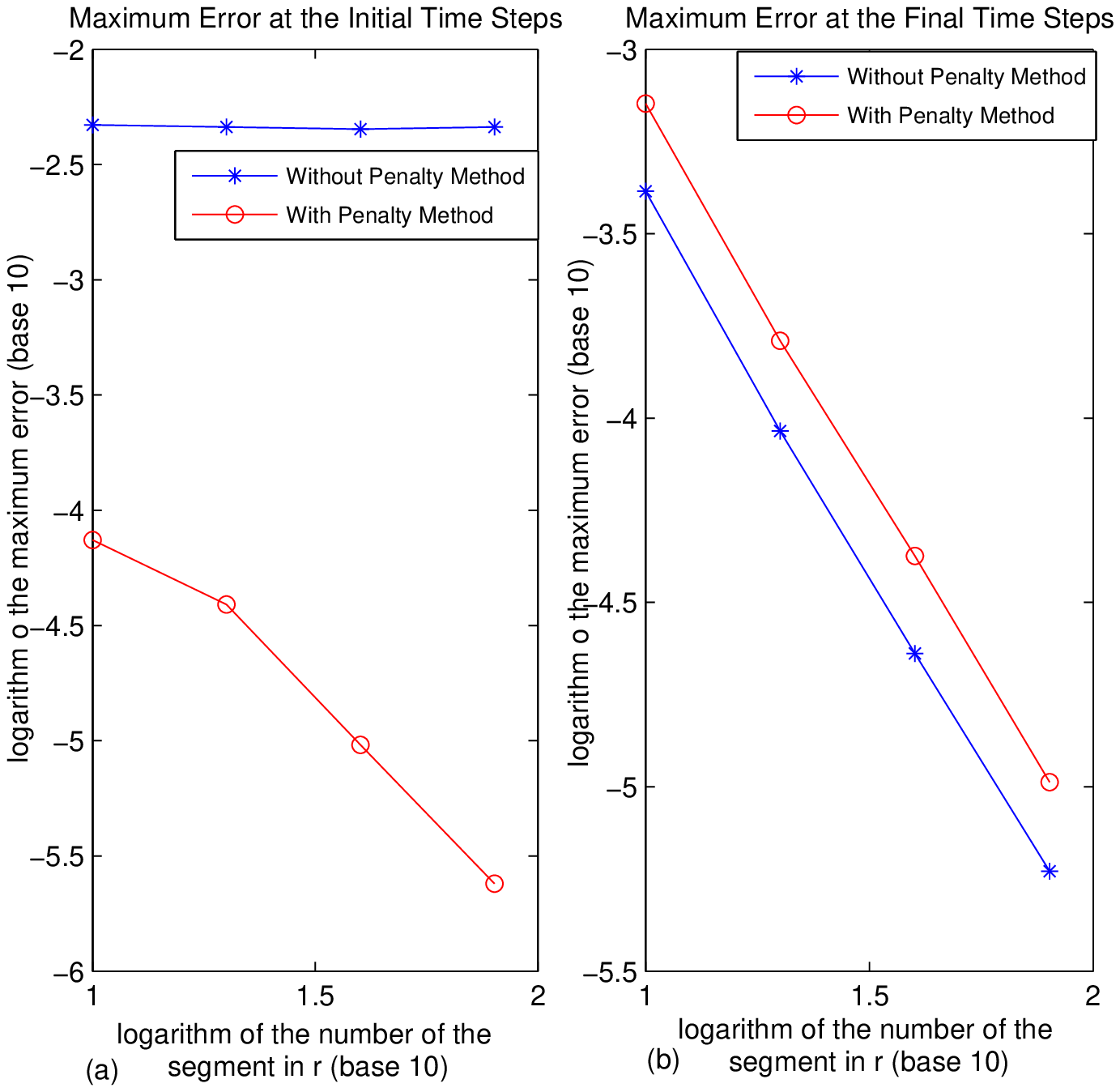}}\caption{Decay of the maximum errors. When we apply the penalty method here $\ep=0.1$. (a) at the initial steps, (b) at the final steps}\label{f12.1}
\end{figure}
Consider the $2D$ system (\ref{e3.5}), where $0\leq t\leq 1$, and set $\nu=0.2,f=0,g=0$ and $u_{0}(x,y)=xy$, so that $g(0)\neq u_{0}|_{\partial \Omega}$. We also set $\varepsilon=0.1$ the same as before. In this case, we face the singularities almost everywhere along the unit circle except at the points where $x=0$ or $y=0$. The effectiveness of this method will be verified with the following numerical results.\\

We first compute the solution of (\ref{e3.5}) without applying the penalty method.  The solution is plotted in
Fig. \ref{f10}; (a) is the graph of the solution at $t=0$, (b) is
the graph of the solution at $t=0.5$, and (c) is the graph of the
solution at $t=1$. As we did for the system (\ref{e3.1}) for the square,
we plot in Fig. \ref{f13}
the sections ($r\in(0,1),\theta=\frac{\pi}{4}$) of the solution
at times close to 0. It is clear that, at $t=0$, the graph displays a sharp gradient around the corner of
the time--space axis due to the
discrepancy between
the initial and boundary conditions there, and as time evolves, the gradient becomes smoother and smoother.\\

To study the error of the system in the disk $\Omega$, we define the maximum comparative errors as for the square $\Omega$. Hence we plot the $L^\infty$ errors for the $2D$ system for the disk $\Omega$ on Fig. \ref{f11}; graph (b) is the maximum comparative error along the whole time period if we apply the penalty method. Because the discrepancy happens at the time-space corner, we zoom into the left corner of graph (b) and compare it with the error when we do not apply the penalty method. From graph (a), we observe that
the magnitude of the errors at the time-space corner are reduced by a factor of $10$ if we apply the penalty method.\\

In Fig.\ref{f12}, we plot the maximum comparative error for (\ref{e3.5}) with a fixed mesh at both initial and final steps. At the initial step, the error decreases sharply as $\varepsilon$ increases and remains close to 0, and then it becomes flat. At the final step, the error increases almost linearly as $\epsilon$ increases. The observation also leads to the following conclusion: at about $\varepsilon=0.1$, the initial error is minimized while the error at final step is well controlled.\\

As for the previous example, we shall now look at how the singularity, induced by the compatibility between the initial and boundary data, affects the convergence rates of the numerical scheme. In Fig. \ref{f12.1} we plot the maximum errors, at the initial and final time steps, with and without the penalty method, against the spatial resolution in the log–log scale. We see in Fig. \ref{f12.1} (a) that, without the penalty method, the maximum errors do not decrease as the grid refines, which demonstrates that the singularity in the solution at the initial time step is serious. With the penalty method, the maximum errors decay at roughly the second order. Fig. \ref{f12.1} (b) shows that, with and without applying penalty method, the maximum errors at the final steps (t=1) decay at approximately the second order.\\

\subsection{Implementation in a 1D System}\label{s3.3}
As we said in the Introduction the penalty method applies without any restriction on space dimension.  However a number of methods have previously been proposed which only apply to space dimension one.  Our aim is now to compare the efficiency of the penalty method with some of the earlier methods; and therefore we can only consider the case of space dimension 1.  More precisely we will consider the Corrector Methods as proposed in \cite{FlFo03}-\cite{FlSw02} and compare them with the penalty method for the $1D$ system
\begin{figure}[ht]
\scalebox{0.5}{\includegraphics{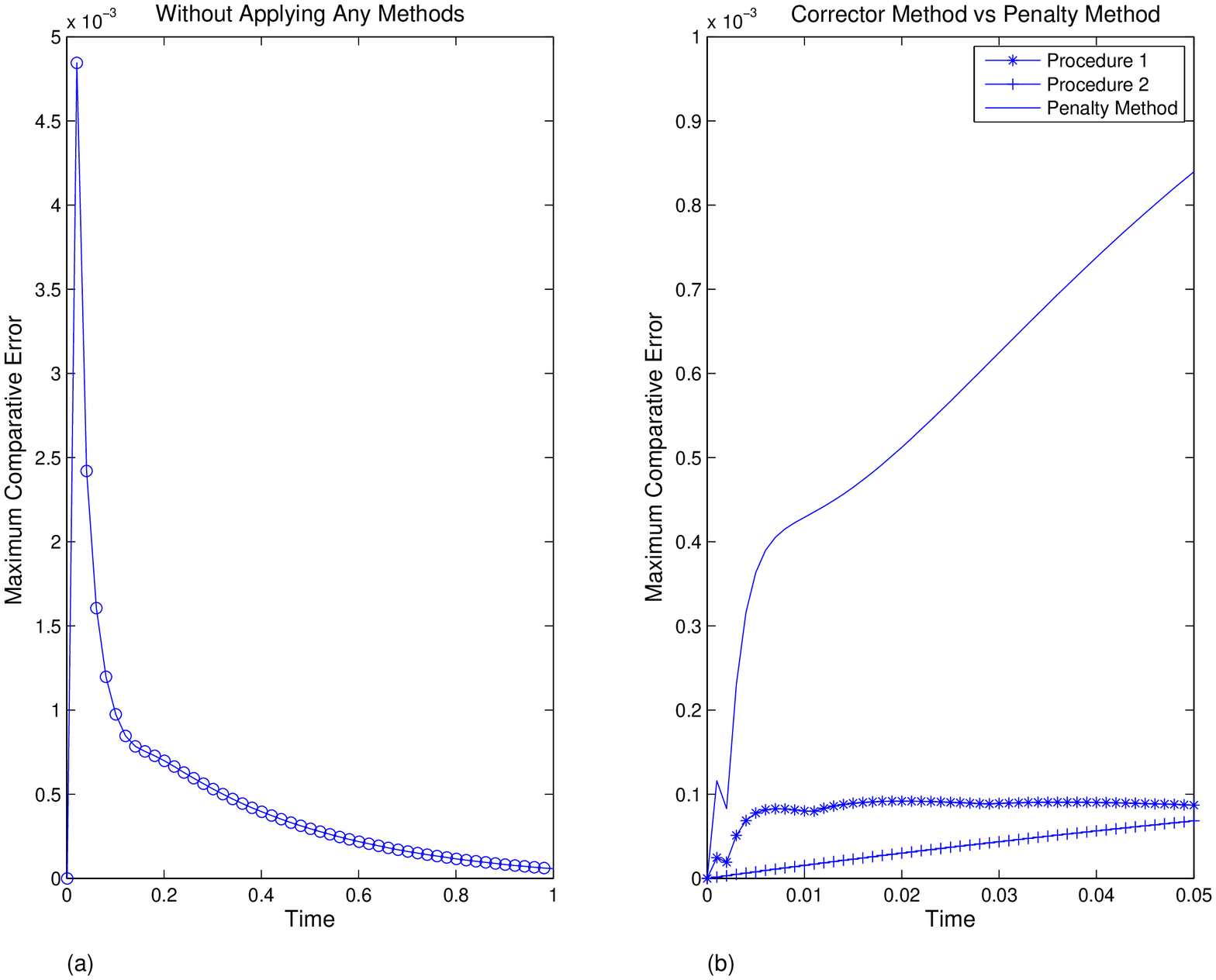}} \caption{Comparative
error of the two methods in 1D system in $L^{\infty}$ sense at
$\varepsilon=0.1$} \label{f3}
\end{figure}
\begin{figure}[ht]
\scalebox{0.55}{\includegraphics{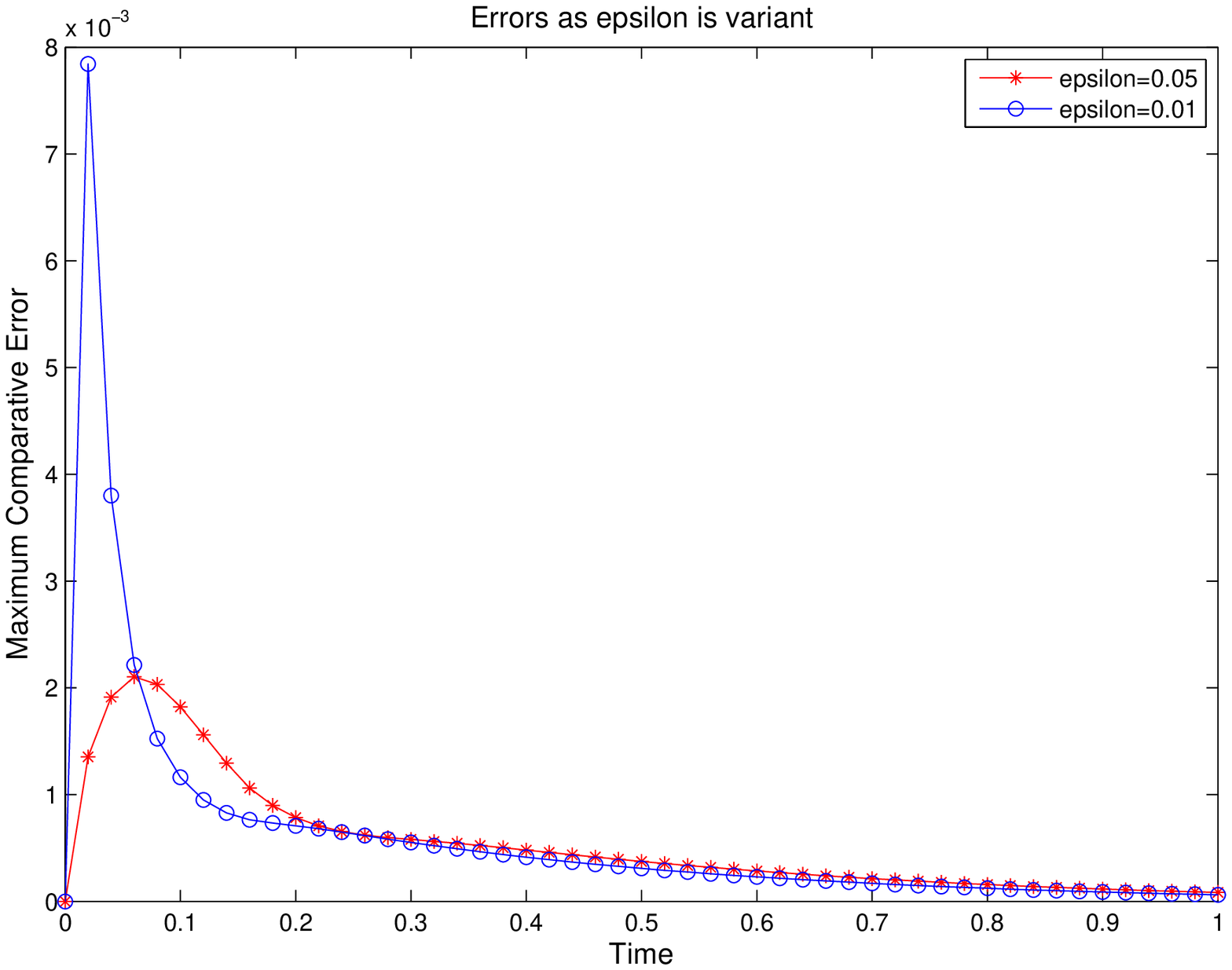}}\caption{The maximum comparative errors for the 1D system in $L^{\infty}$ sense along the time}\label{f7}
\end{figure}

\begin{equation}\left\{\label{e3.2}
\begin{aligned}
&u_{t}-\nu u_{xx}=0, \qquad 0<x<1,\quad 0<t<1,\\
&u(x,0)=u_{0}\\
&u(0,t)=g_{1}(t), \quad u(1,t)=g_{2}(t).
\end{aligned}\right.
\end{equation}
Here we set $u_{0}(x)=sin(\dfrac{5\pi}{4}x+\dfrac{3\pi}{4})$,
$g_{1}(t)=0,\enspace g_{2}(t)=0,\enspace \nu=0.2$. For the Penalty
Method, we also set $\varepsilon=0.1$, and for the Corrector Method, we
have the following choice of correctors \cite{CQT11}-\cite{FlFo04} offering increasing accuracy:
\begin{equation}\label{e3.3}
S=\left\{
\begin{aligned}
&0,\\
&\alpha_{0}S_{0},\hskip0.7in\text{(Procedure 1)}\\
&\alpha_{0}S_{0}+\alpha_{1}S_{1},\hskip0.3in\text{(Procedure 2)}
\end{aligned}\right.
\end{equation}
where $\alpha_{0}=g_{1}(0)-u_{0}(0),\enspace
\alpha_{1}=g_{1t}(0)-u_{0xx}(0)$, $\displaystyle S_{0}=\dfrac{1}{\sqrt{\pi\nu
t}}\int^{\infty}_{x}e^{-\frac{s^{2}}{4\nu t}}ds =
erfc(\dfrac{x}{\sqrt{\nu t}})$ and
$\displaystyle S_{1}=\int^{t}_{0}S_{0}(x,\tau)d\tau$. Here Procedure 1 absorbs the $0^{th}$
order incompatibility ($g_{1}(0)\neq u_{0}(0)$), and Procedure 2 absorbs
both the $0^{th}$ and $1^{st}$ order incompatibilities
($g_{1}(0)\neq u_{0}(0)$ and $g_{1t}(0)\neq \nu u_{0xx}(0)$).

Let $u=v+S$; we see that $v$ is the solution of the following equation
\begin{equation}\left\{\label{e3.6}
\begin{aligned}
&v_{t}-\nu v_{xx}=0,\enspace 0<x<1,\enspace 0<t<1,\\
&v(x,0)=u_{0}(x),\\
&v(0,t)=g_{1}(t)-S(0,t),\enspace v(1,t)=g_{2}(t)-S(1,t).\\
\end{aligned}\right.
\end{equation}

We choose to solve equation (\ref{e3.6}) by finite differences. Fig. \ref{f3} gives the comparison between different methods (Penalty Method and Correction Method).
Figure \ref{f3} (a) gives us the maximum comparative error of system (\ref{e3.2}) without applying any methods. Figure \ref{f3}
(b) compares the two methods, zooming into the corner of the time-space domain where errors are the largest due to the incompatibility at $t=0$. As expected Procedure 2 gives slightly better results than Procedure 1. Also the errors with the penalty method are larger than with both procedures, but still of comparable magnitude whereas the errors without any procedure reach a pick about 6 times larger ($4.8\times 10^{-3}$ vs $0.8\times 10^{-3}$). Now we want to vary $\varepsilon$ in this 1D system, Fig.\ref{f7} shows that if $\varepsilon$ is too small as compared to the mesh, the Penalty Method would not reduce the errors at the spatio-temporal corner, but if it is an appropriate small number, it could really reduce the errors by more than $80\%$.

\section{Conclusion}

The penalty method gives a way to solve the higher dimensional incompatibility problems. As expected, there exists a solution for system (\ref{e1.2}) which is continuous over $[t_{0}, T]$, for all $t_{0}>0$.\\

The discrepancy occurs at the time-space corner; we are effectively interested
in the errors for the initial short time period. The numerical simulations for the system
with both a square $\Omega$ and a disk $\Omega$ yield similar results. At the spatio-temporal corner, the magnitudes of the errors are reduced by about one order of magnitude by the penalty method. Tests are also conducted to study the effects of different values of $\varepsilon$, the key parameter in the penalty method. We find that with an appropriate small value for $\varepsilon$, the initial error can be minimized while the error at final step is under well controlled.\\

Finally, in space dimension one, when both methods are available (penalty method and correction procedures 1 and 2), the penalty method gives a slightly larger error than the Correction Procedures 1 and 2; but the order of magnitude of the errors are comparable and they are all significantly smaller than the errors appearing when no correction procedure is implemented . \\

\section*{Acknowledgments}

This work was partially supported by the National Science Foundation
under the grants NSF-DMS-0604235, and DMS-0906440 and by the
Research Fund of Indiana University.

\vspace{0.5cm}
\begin{center}
{\bf Appendix: The user guide}
\end{center}

The aim is to address the incompatibility issue for the multi-dimensional time-dependent linear parabolic equation
\begin{equation}\left\{\label{e4.1}
\begin{aligned}
&u_{t}-\nu \triangle u= f,\qquad x\in\Omega\subset R^{d}, \quad t\in R^{+},\\
&u|_{t=0}=u_{0},\\
&u|_{\partial\Omega}=g.
\end{aligned}\right.
\end{equation}

where $u_{0}|_{\partial\Omega}\neq g|_{t=0}$. So we consider new system instead, namely, for $\ep>0$ fixed,

\begin{equation}\left\{\label{e4.2}
\begin{aligned}
&u^{\varepsilon}_{t}-\nu\triangle u^{\varepsilon}=f, \qquad x\in\Omega\subset R^{d}, \quad t\in R^{+},\\
&u^{\varepsilon}|_{t=0}=u_{0},\\
&u^{\varepsilon}|_{\partial\Omega}=k^{\varepsilon}.
\end{aligned}\right.
\end{equation}
\begin{equation}\left\{\label{e4.3}
\begin{aligned}
&k^{\varepsilon}_{t}+\frac{1}{\varepsilon}(k^{\varepsilon}-g)=0,\qquad t\in R^{+},\quad\\
&k^{\varepsilon}(0)=u_{0}|_{\partial\Omega}.
\end{aligned}\right.
\end{equation}

We consider for instance the rectangle $0\leq x\leq 1$, $0\leq y\leq 1$ and $0\leq t\leq 1$. We consider the discretization meshes $\triangle x= 1/M$, $\triangle y= 1/N$ and $\triangle t= 1/T$, where $M,N,T$ are integers. We use an explicit scheme to compute the numerical solution of the original system (\ref{e4.1}) and of the modified system (\ref{e4.2}), (\ref{e4.3}), that is respectively:
\begin{equation}\left\{\label{e4.4}
\begin{aligned}
&\f{u^{n+1}_{i,j}-u^{n}_{i,j}}{\triangle t}-\nu(\frac{u^{n}_{i+1,j}+u^{n}_{i-1,j}-2u^{n}_{i,j}}{\triangle x^{2}}+\frac{u^{n}_{i,j+1}+u^{n}_{i,j-1}-2u^{n}_{i,j}}{\triangle y^{2}})=f^{n}_{i,j},\\
&for\enspace 1\leq i\leq N-1,\enspace 1\leq j\leq M-1,\enspace 1\leq n \leq T,\\
&u^{n}_{i,j}|_{\partial\Omega}=g_{i,j}(n\triangle t)|_{\partial\Omega},\enspace for\enspace i=0,N\enspace or\enspace j=0,M,\\
&u^{0}_{i,j}=u_{0}(i\triangle x, j\triangle y), \enspace for \enspace 0\leq i\leq N,\enspace 0\leq j\leq M.
\end{aligned}\right.
\end{equation}
for (\ref{e4.1}), and , for (\ref{e4.2})-(\ref{e4.3}):
\begin{equation}\left\{\label{e4.5}
\begin{aligned}
&\f{u^{n+1}_{i,j}-u^{n}_{i,j}}{\triangle t}-\nu(\frac{u^{n}_{i+1,j}+u^{n}_{i-1,j}-2u^{n}_{i,j}}{\triangle x^{2}}+\frac{u^{n}_{i,j+1}+u^{n}_{i,j-1}-2u^{n}_{i,j}}{\triangle y^{2}})=f^{n}_{i,j},\\
&for\enspace 1\leq i\leq N-1,\enspace 1\leq j\leq M-1,\enspace 1\leq n \leq T,\\
&u^{n}_{i,j}|_{\partial\Omega}=k^{\ep n}_{i,j}|_{\partial\Omega},\enspace for\enspace i=0,N\enspace or\enspace j=0,M,\\
&u^{0}_{i,j}=u_{0}(i\triangle x, j\triangle y), \enspace for \enspace 0\leq i\leq N,\enspace 0\leq j\leq M,\\
&\frac{k^{\ep n+1}_{i,j}-k^{\ep n}_{i,j}}{\triangle t}+\frac{1}{\ep}(k^{\ep n}_{i,j}-g_{i,j}(n\triangle t))=0,\enspace for \enspace i=0,N\enspace or\enspace j=0,M,\enspace n\geq 1,\\
&k^{\ep 0}_{i,j}=u_{0}(i\triangle x, j\triangle y),\enspace for \enspace i=0,N\enspace or\enspace j=0,M.
\end{aligned}\right.
\end{equation}

\bibliographystyle{amsplain}
\bibliography{biblio}

\end{document}